\documentclass{article}

\usepackage{setspace}
  
 \usepackage{hyperref,amssymb,graphicx,epsfig,amsmath,color,caption,subcaption,epstopdf,amsthm}

\newtheorem{definition}{Definition}[section]
\newtheorem{lemma}[definition]{Lemma}
\newtheorem{remark}[definition]{Remark}

\newtheorem{theorem}[definition]{Theorem}

\begin{document}

\title{Feedback Integrators}
\author{Dong Eui Chang \footnote{Corresponding author. Department of Applied Mathematics, University of Waterloo, Waterloo, Ontario N2L 3G1, Canada. dechang@uwaterloo.ca}
 \and 
Fernando Jim\'enez\footnote{Department of Applied Mathematics, University of Waterloo, Waterloo, Ontario, N2L 3G1, Canada. fjimenez@uwaterloo.ca} 
 \and 
 Matthew Perlmutter\footnote{Departamento de Matem\'atica, Universidade Federal de Minas Gerais,
Belo Horizonte, Brazil. matthew@mat.ufmg.br  }
}
\date{\small To appear in {\it Journal of Nonlinear Science}.  \url{http://dx.doi.org/10.1007/s00332-016-9316-7}}
%
%
%
%

%
%

\maketitle

\markright{\hfill FEEDBACK INTEGRATORS\hfill}

\begin{abstract}
A new method is proposed to numerically integrate a dynamical system on a manifold such that the trajectory stably remains on the manifold and preserves first integrals of the system. The idea is that given an initial point in the manifold we extend the dynamics  from the manifold to its ambient Euclidean space and then modify the dynamics outside  the intersection of the manifold and the level sets of the first integrals containing the initial point  such that the intersection becomes a unique local attractor of the resultant dynamics. While the modified dynamics theoretically produces the same trajectory as the original dynamics,  it yields a numerical trajectory  that stably remains on the manifold and preserves the first integrals. The big merit of our method is that the modified dynamics can be integrated with any ordinary numerical integrator such as Euler or Runge-Kutta. We illustrate this method by applying it to three famous problems: the free rigid body,  the Kepler problem and a perturbed Kepler problem with rotational symmetry. We also  carry out simulation studies  to demonstrate the excellence of our method  and make comparisons with the standard projection method, a  splitting method and St\"ormer-Verlet schemes.

\end{abstract}
   
\tableofcontents

\section{Introduction}\label{Intro}
Given a dynamical system on a manifold with first integrals, it is important for a numerical integrator
to preserve the manifold structure and the first integrals of the equations of motion. This has been the focus of much effort in the development of numerical integration schemes \cite{Hair}.
In this paper we do not propose any specific numerical integration scheme, but rather propose a new paradigm of integration that 
can  faithfully preserve  conserved quantities with  existing numerical integration schemes. 

The main idea in our paradigm is as follows.  Consider a dynamical system  on a manifold $M$ with first integrals $f_i : M \rightarrow \mathbb R$, $i = 1, \ldots, \ell$. Assume that we can embed the manifold $M$ into Euclidean space $\mathbb R^n$ and extend the first integrals to a neighborhood $U$ of $M$ in $\mathbb R^n$.  For an arbitrary point $x_0 \in M$, consider the set
\[
\Lambda = \{ x \in U \mid x \in M, f_i(x) = f_i(x_0), i = 1, \ldots, \ell\}
\]
which is  the intersection of $M$ with all the level sets of the first integrals containing the point $x_0$, and is an invariant set of the dynamical system. We then extend the dynamical system  from $M$ to $U$ and then modify the dynamics outside of $\Lambda$
such that the set $\Lambda$
becomes a unique local attractor of the extended, modified system. Since the dynamics have not changed on $\Lambda$ by the extension and modification to $U$, both the original system on $M$ and the extended, modified system on $U$  produce the same trajectory for the initial point $x_0 \in \Lambda$. Numerically, however, integrating the extended system has the following advantage: if the trajectory deviates from $\Lambda$ at some numerical integration step, then it will get pushed back toward the attractor $\Lambda$ in the extended, modified dynamics, thus remaining on the manifold $M$ and preserving all the first integrals.   It can be rigorously shown that the discrete-time dynamical system derived from any one-step numerical integrator with uniform step size $h$ for the extended, modified  continuous-time system  indeed has an attractor $\Lambda_h$ that contains the set $\Lambda$ in its interior and converges to $\Lambda$ as $h \rightarrow 0+$. In this paper we shall use the word, {\it preserve}, in this sense.
It is noteworthy that the numerical integration of the extended dynamics 
can be carried out  with any ordinary  integrator and is done in one global Cartesian coordinate system on $\mathbb R^n$. We  find  conditions for applicability of this method and implement the result  on the following three examples: the free rigid body dynamics, the Kepler problem, and a  perturbed Kepler problem with rotational symmetry.  We  also carry out simulation studies  to show the excellence of our new paradigm of integration for numerical preservation of conserved quantities in comparison with other well-known integration schemes, such as projection and splitting methods and symplectic St\"ormer-Verlet integrators.

\section{Theory}\label{section:theory}

Consider a dynamical system on an open subset $U$ of $\mathbb R^n$:
\begin{equation}\label{our:dyn:sys}
\dot x = X(x),
\end{equation}
where $X$ is a $C^1$ vector field on $U$. Let us make the following assumptions:
\begin{itemize}
\item [A1.]  There is a $C^2$ function $V : U \rightarrow \mathbb R$ such that  $V(x) \geq 0$ for all $x\in U$, $V^{-1}(0) \neq \emptyset$, and 
\begin{equation}\label{Vdot:X}
\nabla V(x) \cdot X (x) = 0
\end{equation}
for all $x \in U$.
 
 \item [A2.] There is a positive number $c$ such that $V^{-1}([0,c])$ is a compact subset of $U$.
 
\item [A3.] The set of  all critical points of $V$  in $V^{-1}([0,c])$ is equal to $V^{-1}(0)$. 
\end{itemize}
Adding the negative gradient of $V$ to (\ref{our:dyn:sys}), let us consider the following dynamical system on $U$:
\begin{equation}\label{our:new:dyn:sys}
\dot x = X(x) - \nabla V(x).
\end{equation}
Since $0$ is the minimum value of $V$, $\nabla V (x)= 0$ for all $x\in V^{-1}(0)$. Hence, the two vector fields $X$ and $X-\nabla V$ coincide on $V^{-1}(0)$.
\begin{theorem}\label{theorem:general}
Under assumptions A1 -- A3, every trajectory of (\ref{our:new:dyn:sys}) starting from a point in $V^{-1}([0,c])$ stays in $V^{-1}([0,c])$ for all $t\geq 0$ and asymptotically converges to the  set $V^{-1}(0)$ as $t \rightarrow \infty$. Furthermore, $V^{-1}(0)$ is an invariant set of   both  (\ref{our:dyn:sys}) and (\ref{our:new:dyn:sys}).
\begin{proof}
  Let $x(t)$ be a trajectory of (\ref{our:new:dyn:sys}) starting from a point in $V^{-1}([0,c])$. By A1
\begin{equation}\label{V:dot}
\frac{d}{dt} V(x(t)) = \nabla V (x(t)) \cdot ( X(x(t))- \nabla V(x(t)))= -|\nabla V(x)|^2 \leq 0
\end{equation}
for all $t $. Hence, $V^{-1}([0,c])$ is a positively invariant set of (\ref{our:new:dyn:sys}). From  (\ref{V:dot}) and A3, it follows that $\{  x \in V^{-1}([0,c]) \mid \dot V (x) = 0\} = \{x \in V^{-1}([0,c]) \mid \nabla V(x) = 0\}=  V^{-1}(0)$. Hence, by LaSalle's invariance principle \cite{La60}, $x(t)$ converges asymptotically to $V^{-1}(0)$ as $t \rightarrow \infty$, where A2 is used for compactness of $V^{-1}([0,c])$.  The invariance of $V^{-1}(0)$ follows from (\ref{Vdot:X}) and the coincidence of (\ref{our:dyn:sys}) and (\ref{our:new:dyn:sys}) on $V^{-1}(0)$.
\end{proof}
\end{theorem}

Let us find a higher-order condition than that in assumption A3 so that A3 can be relaxed. 
For the function $V$  and the vector field $X$ in the statement of assumption A1, which are  now both assumed to be of $C^\infty$, let 
\begin{equation}\label{def:set:S}
S = \left \{ x \in U  \bigg{|} X^k \frac{\partial V}{\partial x^i} = 0 \,\, \forall\,\, k \geq 0, 1 \leq i \leq n  \right \},  
\end{equation}
where $x = (x^1, x^2, \ldots, x^n)$, and $X^k \frac{\partial V}{\partial x^i} $ denotes the $k-$th order directional derivative of $\partial V/ \partial x^i$ along $X$, i.e.,
\[
X^0 \frac{\partial V}{\partial x^i} = \frac{\partial V}{\partial x^i}; \quad X \frac{\partial V}{\partial x^i}  =X \cdot  \nabla \frac{\partial V}{\partial x^i}; \quad  X^k \frac{\partial V}{\partial x^i} =  X\left (X^{k-1} \frac{\partial V}{\partial x^i} \right), \,k \geq 2.
\]
Consider  the following assumption in place of A3:
\begin{itemize}
\item [A3$^{\prime}$.] $S \cap V^{-1}([0,c]) \subset V^{-1}(0)$.
\end{itemize}
The following theorem generalizes Theorem \ref{theorem:general}:
\begin{theorem}\label{theorem:general:more}
Under assumptions A1, A2 and  A3$^{\prime}$, every trajectory of (\ref{our:new:dyn:sys}) starting  in $V^{-1}([0,c])$  stays in $V^{-1}([0,c])$ for all $t\geq 0$ and asymptotically converges to the  set $V^{-1}(0)$ as $t \rightarrow \infty$. Furthermore, $V^{-1}(0)$ is an invariant set of   both  (\ref{our:dyn:sys}) and (\ref{our:new:dyn:sys}).
\begin{proof}
Consider the dynamics (\ref{our:new:dyn:sys}). It is easy to show that $V^{-1}([0,c])$ is a positively invariant set of the dynamics.
Let $\mathcal M$ be the largest invariant set  in $\mathcal E = \{x \in U\mid  \dot V (x)= 0\} \cap V^{-1}([0,c])$.  Let $x(t)$ be an arbitrary trajectory in $\mathcal M$.  Since $\mathcal E= \{x \in U \mid  \nabla V (x)= 0\} \cap V^{-1}([0,c])$ as shown  in the proof of Theorem \ref{theorem:general},   the trajectory $x(t)$  satisfies $\nabla V = 0$, i.e.,
\begin{equation}\label{dVdxi}
\frac{\partial V}{\partial x^i} (x(t)) =0
\end{equation}
for all $t \in \mathbb R$ and $ 1\leq i \leq n$. Since $\nabla V = 0$ along $x(t)$, the trajectory $x(t)$  satisfies
\begin{equation}\label{reduced:X}
\dot x(t) = X(x(t))
\end{equation}
for all $t \in \mathbb R$.
By  differentiating (\ref{dVdxi}) repeatedly in $t$ and using  (\ref{reduced:X}) on each differentiation, we can show that the trajectory $x(t)$ satisfies
\[
X^k \frac{\partial V}{\partial x^i} = 0
\]
for all $t\in \mathbb R$,  $k \geq 0$ and $ 1\leq i \leq n$. Thus, the entire trajectory $x(t)$ is contained in the set $S$ defined in (\ref{def:set:S}), implying $\mathcal M \subset S$, from which and  A3$^{\prime}$  it follows $\mathcal M \subset V^{-1}(0)$. Hence, by LaSalle's invariance principle, every trajectory starting in $V^{-1}([0,c])$ asymptotically converges to $\mathcal M$ and thus to $V^{-1}(0)$ as $t \rightarrow \infty$.

The invariance of $V^{-1}(0)$ follows from (\ref{Vdot:X}) and the coincidence of (\ref{our:dyn:sys}) and (\ref{our:new:dyn:sys}) on $V^{-1}(0)$.\end{proof}
\end{theorem}

\begin{remark}
1. If condition (\ref{Vdot:X}) is replaced by $\nabla V(x)\cdot X(x) \leq 0$ in assumption A1, then Theorems  \ref{theorem:general} and \ref{theorem:general:more}  still hold provided that the invariance of $V^{-1}(0)$ is replaced by positive invariance in the statement of the theorems.

2.Theorems \ref{theorem:general} and \ref{theorem:general:more} still hold with 
the use of the following modified dynamics
\[
\dot x = X(x) - A(x) \nabla V(x)
\]
 instead of (\ref{our:new:dyn:sys}), where $A(x)$ is an $n\times n$ matrix-valued function with its symmetric part $(A(x) + A^T(x))$  positive definite at each $x\in\mathbb R^n$.

3. From the control viewpoint, the added term $-\nabla V(x)$ in (\ref{our:new:dyn:sys}) can be regarded as a negative feedback control $ u (x) =-\nabla V(x)$ to asymptotically stabilize the set $V^{-1}(0)$ for the control system $\dot x = X(x) + u$ with control $u$.
\end{remark}

Suppose that assumptions A1, A2 and  A3 (or  A3$^{\prime}$ instead of A3) hold and that we want to integrate the dynamics  (\ref{our:dyn:sys}) for an initial point $x(0) \in V^{-1}(0)$. Since $V^{-1}(0)$ is positively invariant, the trajectory must remain in $V^{-1}(0)$ for all $t\geq 0$. Recall that the two dynamics (\ref{our:dyn:sys})  and (\ref{our:new:dyn:sys}) coincide on $V^{-1}(0)$, so  we can integrate (\ref{our:new:dyn:sys}) instead of (\ref{our:dyn:sys}) for the initial condition. Though there is no theoretical difference between the two integrations, integrating  (\ref{our:new:dyn:sys})  has a numerical advantage over integrating (\ref{our:dyn:sys}). Suppose that the trajectory numerically deviates from the positively invariant set $V^{-1}(0)$ during integration.  Then the dynamics  (\ref{our:new:dyn:sys}) will push the trajectory back toward $V^{-1}(0)$ since $V^{-1}(0)$ is the attractor of  (\ref{our:new:dyn:sys})  in  $V^{-1}([0,c])$  whereas the dynamics  (\ref{our:dyn:sys}) will  leave the trajectory outside of $V^{-1}(0)$ which would  not happen  in the exact solution. It is  noteworthy that this integration strategy is independent of the choice of  integration schemes.   In the Appendix we show  that  any one-step numerical integrator, as a discrete-time dynamical system, with uniform step size $h$  for (\ref{our:new:dyn:sys}) has an attractor $\Lambda_h$  that contains $V^{-1}(0)$ in its interior and converges to $V^{-1}(0)$ as $h\rightarrow 0+$.

Let us  now apply this integration strategy to numerically integrate dynamics on a manifold while preserving its first integrals and the domain manifold. 
Consider a manifold $M$ and dynamics 
\begin{equation}\label{dyn:on:M}
\dot x = X(x)
\end{equation}
on $M$ that have $\ell$ first integrals $f_i :M \rightarrow \mathbb R$, $i = 1, \ldots, \ell$.    Suppose that $M$ is an embedded manifold in $\mathbb R^n$ as a  level set of a function $f_0: \mathbb R^n \rightarrow \mathbb R^r$ for some $r$, and that both the dynamics (\ref{dyn:on:M}) and the functions $f_i$, $i = 0, \ldots, \ell$ extend to an open neighborhood $U$ of $M$ in $\mathbb R^n$. Our goal is to numerically integrate (\ref{dyn:on:M}) with an initial condition $x(0) = x_0 \in M$ while preserving the manifold $M$ and the first integrals. Let 
\begin{equation}\label{def:f}
f = (f_0, f_1, \ldots, f_\ell): \mathbb R^n \rightarrow \mathbb R^{r+\ell}
\end{equation}
 and define a function $V:U \subset \mathbb R^n \rightarrow \mathbb R$ by
\begin{equation}\label{def:V:intro}
V(x) = \frac{1}{2} (f(x) - f(x_0))^T K (f(x) - f(x_0)),
\end{equation}
where $K$ is an $(r+\ell) \times (r+\ell)$ constant symmetric positive definite matrix.   Notice that
\[
V^{-1}(0) = \{x \in U \mid x \in M, f_i(x) = f_i(x_0), i = 1, \ldots, \ell\},
\]
and that $V^{-1}(0)$ is invariant under the flow of (\ref{dyn:on:M}). Or, more generally we can define  a function $V(x)$ as $V(x) = W(f_0(x), f_1(x), \ldots, f_\ell(x))$ where $W: \mathbb R^{r +\ell } \rightarrow \mathbb R$ is a non-negative  function that takes the  value of $0$ only at $(f_0(x_0), f_1(x_0), \ldots, f_\ell(x_0))$. 
If  the function $V$ satisfies assumptions A1, A2 and  A3 (or  A3$^{\prime}$ instead of A3), then by Theorem \ref{theorem:general} (or Theorem \ref{theorem:general:more}), $V^{-1}(0)$ is the local attractor of the modified dynamics
\begin{equation}\label{dyn:on:M:ext}
\dot x = X(x) -\nabla V(x)
\end{equation}
that coincide with the original dynamics (\ref{dyn:on:M})  on $V^{-1}(0)$. 

The following lemma  provides a sufficient condition under which  the function $V$ defined in (\ref{def:V:intro})  satisfies assumptions A2 and A3:

\begin{lemma}\label{kuhl:lemma}
Consider the functions $f$ and $V$ defined in (\ref{def:f}) and (\ref{def:V:intro}). If  $V^{-1}(0)$ is compact and the Jacobian matrix $Df(x)$ of $f$ has rank $(r+\ell)$ for all  $x \in V^{-1}(0)$, then there is a number $c>0$ such that assumptions A2 and A3 hold. 
\begin{proof}
By compactness of $V^{-1}(0)$ and the regularity of $Df$, there is a bounded open set $X$ such that  $V^{-1}(0) \subset X \subset \operatorname{cl}( X) \subset U$, and $Df$(x) has rank $r+\ell$ for all  $x \in X$, where $\operatorname{cl}(X)$ denotes the closure of $X$. Consider now the gradient of $V$. An easy calculation shows that,
\[
\nabla V(x) = Df(x)^TK (f(x) - f(x_0)).
\] 
Now, since for all $x\in X$, $Df(x)$ is onto as a linear map, $Df(x)^{T}$ is therefore one to one. It follows that, for $x\in X$,
\begin{equation}\label{onetooneDfT}
\nabla V(x)=0 \iff  f(x)-f(x_{0})=0 \iff x\in V^{-1}(0).
\end{equation}
In other words, the set of all critical points of $V$ in $X$ is equal to $V^{-1}(0)$. Since the boundary $\partial X$ of $X$, being closed and bounded, is compact and $\partial X \cap V^{-1}(0) = \emptyset$, the minimum value, denoted by $d$, of $V$ on $\partial X$ is positive.  If necessary, restrict the function $V$ to $X$, replacing its original domain $U$ with $X$. Then, there is a positive number $c$ less than $d$ such that $V^{-1}([0, c]) \subset X$. Therefore,  assumption A3 holds for this number $c$. Since the closed set $V^{-1}([0,c])$ is contained in the bounded set $X$, it is compact, which implies that assumption A2 holds. 
\end{proof}
\end{lemma}
\begin{theorem}\label{kuhl:theorem}
For the functions $f$ and $V$ defined in  (\ref{def:f}) and (\ref{def:V:intro}), if  $V$ satisfies (\ref{Vdot:X}) for all $x\in U$, the set $V^{-1}(0)$ is compact and the Jacobian matrix $Df(x)$ is onto  for all  $x \in V^{-1}(0)$, then there is a number $c>0$ such that every trajectory starting in $V^{-1}([0,c])$  remains in $V^{-1}([0,c])$ for all $t\geq 0$ and asymptotically converges to $V^{-1}(0)$ as $t\rightarrow \infty$. 
\end{theorem}

\begin{theorem}
For the functions $f$ and $V$ defined in (\ref{def:f}) and (\ref{def:V:intro}), if  $V$ satisfies (\ref{Vdot:X}) for all $x\in U$, the set $V^{-1}(0)$ is compact and there is an open subset $X$ of $U$ containing $V^{-1}(0)$ such that the Jacobian matrix $Df(x)$ is onto  for all  $x \in X\backslash V^{-1}(0)$, then there is a number $c>0$ such that every trajectory starting in $V^{-1}([0,c])$  remains in $V^{-1}([0,c])$ for all $t\geq 0$ and asymptotically converges to $V^{-1}(0)$ as $t\rightarrow \infty$. 
\begin{proof}
Modify the proof of Lemma~\ref{kuhl:lemma} appropriately.
\end{proof}
\end{theorem}

As discussed above, we can  integrate (\ref{dyn:on:M:ext}) instead of  (\ref{dyn:on:M}) for the initial condition $x(0) = x_0 \in V^{-1}(0)$, which will yield a trajectory that is expected to numerically well remain on  the manifold $M$ and preserve the values of the first integrals $f_i$, $i = 1, \ldots, \ell$.  It is noteworthy  that the integration is carried out in one Cartesian coordinate system on  $\mathbb R^n$ rather than over local charts on the manifold $M$ which would take additional computational costs for coordinate changes between local charts.  
In the following section, we will apply this strategy  to  the free rigid body dynamics, the Kepler problem and a perturbed Kepler problem with rotational symmetry to integrate the dynamics preserving their first integrals and domain manifolds.

\section{Applications}

\subsection{The Free Rigid Body} 
\label{subsection:rigid:body}

Consider the free rigid body dynamics:
\begin{subequations}\label{RigidBody}
\begin{align}
\dot R &= R\,\hat \Omega, \label{free:rigid:R} \\ 
\dot \Omega &=  {\mathbb I}^{-1} \left((\mathbb I \Omega) \times \Omega \right), \label{free:rigid:Omega}
\end{align}
\end{subequations}  
where $(R, \Omega) \in {\operatorname{SO}(3)} \times \mathbb R^3$; $\mathbb I$ is the moment of inertia matrix;  and 
\begin{equation}\label{def:hat:map}
\hat \Omega = \begin{bmatrix}
0 & -\Omega_3 & \Omega_2 \\
\Omega_3 &0 & -\Omega_1 \\
-\Omega_2 & \Omega_1 & 0
\end{bmatrix}
\end{equation}
for 
\[
\Omega = \begin{bmatrix} \Omega_1 \\  \Omega_2 \\ \Omega_3 \end{bmatrix}.
\] 
Since ${\operatorname{SO}(3)} \subset \mathbb R^{3\times 3}$, from here on we assume that the rigid body dynamics are defined on the Euclidean space $\mathbb R^{3\times 3} \times \mathbb R^3$ and that the matrix $R$ denotes a  $3\times 3$ matrix, not necessarily in ${\operatorname{SO}(3)}$.  This is the {\em extension} of the dynamics step.

  Define  two functions $E: \mathbb R^3 \rightarrow \mathbb R$ and  $\pi: \mathbb R^{3\times 3} \times \mathbb R^3 \rightarrow \mathbb R^3$  by 
\begin{align}
E(\Omega) &= \frac{1}{2} \Omega^T\mathbb I \Omega, \label{def:E}\\
\pi (R, \Omega) &= R\,\mathbb I\, \Omega, \label{def:pi}
\end{align}
where $E$ represents the kinetic energy of the free rigid body and $\pi$ the spatial angular momentum vector when $R \in \operatorname{SO}(3)$.  These quantities are first integrals of (\ref{RigidBody}).  Choose any
\[
R_0 \in \operatorname{SO}(3), \quad \Omega_0 \in \mathbb R^3\backslash \{(0,0,0)\},
\]
and let
\begin{equation}\label{def:E0:pi0}
E_0 = E(\Omega_0) >0, \quad \mathbb \pi_0 = \pi(R_0,\Omega_0) \in \mathbb R^3\backslash \{(0,0,0)\}.
\end{equation}
Define an open set $U$ by
\[
U = \{ (R,\Omega) \in \mathbb R^{3\times3} \times \mathbb R^3 \mid \det (R) >0\}
\]
and a function $V : U \subset \mathbb R^{3\times 3} \times \mathbb R^3 \rightarrow \mathbb R$ by 
\begin{equation}\label{def:V:rigid}
V(R,\Omega) = \frac{k_0}{4}\| R^TR - I\|^2 +  \frac{k_1}{2} \left |  E(\Omega) - E_0 \right | ^2+ \frac{k_2}{2}|\pi (R,\Omega) - \pi_0|^2
\end{equation}
for $(R,\Omega) \in U \subset \mathbb R^{3\times 3} \times \mathbb R^3$, where $k_i >0$, $i = 0, 1, 2$ are constants, and  $\| \cdot \|$ is the 2-norm defined by $\| A\| = \sqrt{\operatorname{trace}(A^TA)}$ for a matrix $A$.  Observe that we are endowing the space $\mathbb{R}^{3\times3}\times\mathbb{R}^{3}$ with the standard inner product, and that the trace norm is precisely the norm induced on $\mathbb{R}^{3\times3}$ by this inner product. We compute all gradients that follow with respect to this inner product.   Notice that
\[
V^{-1}(0) = \{ (R,\Omega) \in \mathbb R^{3\times 3} \times \mathbb R^3 \mid R \in \operatorname{SO}(3), \, E(\Omega) = E_0, \, \pi(R,\Omega) = \pi_0\}. 
\]
\begin{lemma}\label{lemma:gradient}
The gradient $(\nabla_R V, \nabla_\Omega V) \in \mathbb R^{3\times 3} \times \mathbb R^3$ of the function $V$ (\ref{def:V:rigid}) is given by
\begin{subequations}\label{GradEquations}
\begin{align}
\nabla_R V&= k_0R(R^TR - I) + k_2 (\pi(R,\Omega) - \pi_0) \Omega^T \mathbb I, \label{GradEquationsa}\\
\nabla_\Omega V  &=k_1 (E(\Omega)-E_0) \mathbb I \Omega + k_2\mathbb IR^T (\pi(R,\Omega) - \pi_0).\label{GradEquationsb}
\end{align}
\end{subequations}
\begin{proof}
Straightforward.
\end{proof}
\end{lemma}

The following lemma shows that the function $V$  satisfies assumption A1 stated in \S\ref{section:theory}.
\begin{lemma}\label{lemma:1:rigid}
The function $V$ satisfies 
\begin{equation}\label{ok:chef}
\langle (\nabla_R V, \nabla_\Omega V), (R\hat \Omega, \mathbb I ^{-1} ( (\mathbb I \Omega) \times \Omega ) )\rangle = 0.
\end{equation}
\begin{proof}
 One can compute
\begin{align*}
\langle \nabla_R V,  (R\hat \Omega) \rangle&= \operatorname{trace} (\hat \Omega^T R^T (k_0R(R^TR - I) + k_2 (\pi(R,\Omega) - \pi_0) \Omega^T \mathbb I) ) \\
&= -k_0\operatorname{trace} (\hat \Omega R^T R(R^TR - I) ) - k_2\operatorname{trace} (\hat \Omega R^T (\pi(R,\Omega) - \pi_0) \Omega^T \mathbb I )\\
&= -k_2  \Omega^T \mathbb I \hat \Omega R^T (\pi(R,\Omega) - \pi_0),
\end{align*}
where, in the third equality we use the fact that for $A$ symmetric and $B$ antisymmetric, $\operatorname{trace}(AB)=0$.

Next, we compute,
\begin{align*}
\langle \nabla_\Omega V , \mathbb I ^{-1} ( (\mathbb I \Omega) \times \Omega ) \rangle  &= \langle k_1 (E(\Omega)-E_0) \mathbb I \Omega + k_2\mathbb IR^T (\pi(R,\Omega) - \pi_0),   \mathbb I ^{-1} ( (\mathbb I \Omega) \times \Omega ) \rangle \\
&=k_1 (E(\Omega)-E_0) \langle (\mathbb I \Omega) \times \Omega , \Omega \rangle + k_2 \langle  R^T (\pi(R,\Omega) - \pi_0), (\mathbb I \Omega) \times \Omega \rangle \\
&=  k_2 \langle \mathbb I \Omega, \Omega \times R^T (\pi(R,\Omega) - \pi_0) \rangle\\
&= k_2  \Omega^T \mathbb I \hat \Omega R^T (\pi(R,\Omega) - \pi_0).
\end{align*}
Hence,
\[
\langle (\nabla_R V, \nabla_\Omega V), (R\hat \Omega, \mathbb I ^{-1} ( (\mathbb I \Omega) \times \Omega ) )\rangle =  \langle \nabla_R V,  R\hat \Omega \rangle + \langle  \nabla_\Omega V, \mathbb I ^{-1} ( (\mathbb I \Omega) \times \Omega ) \rangle=0.
\]

\end{proof}
\end{lemma}

The following lemma shows that the function $V$ satisfies assumptions A2 and A3  stated in \S\ref{section:theory}.
\begin{lemma}\label{lemma:2:rigid}
There is a number $c$ satisfying
\begin{equation}\label{inequality:c}
0< c < \min \{ k_0/4, k_1|E_0|/2, k_2|\pi_0|^2/2\}
\end{equation}
 such that   $V^{-1}([0,c])$ is a compact subset of $U$ and the set of all critical points of $V$ in $V^{-1}([0,c])$ is equal to $V^{-1}(0)$. 
 \begin{proof}
 It is obvious that there is a number $c$  satisfying  (\ref{inequality:c}) such that $V^{-1}([0,c])$ becomes a compact set in $U$. For such a number $c$, the matrix $R$ is invertible for every $(R,\Omega) \in V^{-1}([0,c])$.
Since $0$ is the minimum value of $V$, every point in $V^{-1}(0)$ is a critical point of $V$. 

Let $(R,\Omega) $ be a critical point of $V$ in  $V^{-1}([0,c]) \backslash V^{-1}(0)$. By Lemma \ref{lemma:gradient} it  satisfies
\begin{subequations}\label{Critical}
\begin{align}
k_0R(R^TR - I) + k_2 (\pi- \pi_0) \Omega^T \mathbb I &=0,  \label{critical:1a} \\
k_1 (E-E_0) \mathbb I \Omega + k_2\mathbb IR^T (\pi- \pi_0) &= 0, \label{critical:1b}
\end{align}
\end{subequations}
where 
\[
\pi = \pi(R,\Omega), \quad E = E(R,\Omega).
\]
Post-multiplying (\ref{critical:1a}) by $R^T$ and pre-multiplying (\ref{critical:1b}) by $R$ yield
\begin{subequations}\label{Critical2}
\begin{align}
k_0R(R^TR - I) R^T+ k_2 (\pi- \pi_0) \pi^T &=0,  \label{critical:2a} \\
k_1 (E-E_0) \mathbb \pi + k_2R\mathbb IR^T (\pi- \pi_0) &= 0, \label{critical:2b}
\end{align}
\end{subequations}
since $\pi =R\mathbb I \Omega$. 
Notice that  $\Omega= 0$ would imply $V(R,\Omega) \geq \frac{k_2}{2}|\pi_0|^2 >c$, contradicting $(R,\Omega) \in V^{-1}([0,c])$. Hence, $\Omega \neq 0$. It follows  from  (\ref{Critical}) that  if any of the three equations
\[
R^TR - I = 0, \quad \pi - \pi_0 = 0, \quad E - E_0 = 0
\]
holds, then the three of them all hold. Thus
\begin{equation}\label{contra:eqn}
R^TR \neq I, \quad \pi \neq \pi_0, \quad E \neq E_0
\end{equation}
since $(R,\Omega) \notin V^{-1}(0)$.
Since the matrix $(\pi- \pi_0) \Omega^T \mathbb I$ in (\ref{critical:1a}) has rank 1 and the matrix $(R^TR - I) $ is symmetric, there exist a unit vector $u \in \mathbb R^3$ and a  number $\kappa \neq 0$ such that
\begin{equation}\label{RTR_I}
R^TR -I = \kappa u u^T.
\end{equation}
Substitution of (\ref{RTR_I}) into (\ref{critical:1a}) and  (\ref{critical:2a})  yields
\begin{align*}
k_0\kappa Ru u^T+ k_2 (\pi -\pi_0) \Omega^T \mathbb I= 0,\\
k_0\kappa Ruu^TR^T + k_2 (\pi -\pi_0) \pi^T = 0,
\end{align*}
which implies 
\begin{equation}\label{parallel_rels}
Ru \parallel \pi \parallel \pi_0, \quad u \parallel \mathbb I \Omega,
\end{equation}
where the symbol $\parallel$ means `is parallel to.'
Hence, we can express $R$ and $\pi$ as 
\begin{align}
R &=  w_1 u_1^T + w_2u_2^T +  a e_{\pi_0} u^T,  \label{R_split} \\
 \pi &= b \pi_0,\label{pi_b_pi0}
\end{align}
for some numbers $a \neq 0$, $b \neq 1$ and vectors $u_1, u_2, w_1, w_2 \in \mathbb R^3$, where $e_{\pi_0} = \pi_0 / |\pi_0|$ and the vectors $u_1$ and $u_2$ can be any vectors such that $\{u_1, u_2, u\}$ becomes an orthonormal basis for $\mathbb R^3$. Substitution of (\ref{R_split})  into (\ref{RTR_I}) implies that $\{ w_1,w_2, e_{\pi_0}\}$ is an orthonormal basis for  $\mathbb R^3$. Substitution of (\ref{R_split}) and (\ref{pi_b_pi0}) into (\ref{critical:1b}) implies $\mathbb I \Omega \parallel \mathbb I u$, which together with $u \parallel \mathbb I \Omega $ in (\ref{parallel_rels}), implies $u \parallel \mathbb Iu$, i.e., $u$ is an eigenvector of $\mathbb I$.  We can now  choose or re-define the unit vectors $u_1$ and $u_2$ such that they become eigenvectors of the symmetric  matrix $\mathbb I$, too. 
In the orthonormal basis $\{u_1, u_2, u\}$, we can now write the moment of inertia matrix $\mathbb I$ as
\[
\mathbb I = I_1 u_1 u_1 ^T + I_2 u_2 u_2^T + I_3 uu^T,
\]
where $I_1, I_2, I_3$ are the eigenvalues of $\mathbb I$, which are all positive, corresponding to the eigenvectors $u_1, u_2, u$, respectively. It is  then easy to see that equations (\ref{Critical2}) imply
\begin{subequations}\label{abEquations}
\begin{align}
k_0 a^2 (a^2 -1) + k_2 |\pi_0|^2b(b-1) &=0,  \label{abEquations:a} \\
k_1 \left ( \frac{|\pi_0|^2 b^2}{2I_3a^2 } - E_0 \right )b + k_2 I_3 a^2 (b-1) &=0,\label{abEquations:b}
\end{align}
\end{subequations}
where we have used $E = (1/2)\Omega^T\mathbb I \Omega = (1/2) \pi^T (R\mathbb I R^T)^{-1}\pi = |\pi_0|^2b^2/2I_3a^2$.

We consider the following two separate cases: $E_ 0= |\pi_0|^2/2I_3$ and  $E_ 0\neq |\pi_0|^2/2I_3$.
Suppose $E_ 0= |\pi_0|^2/2I_3$.  If $b\leq 0$, then 
\[
V(R,\Omega) \geq \frac{k_2}{2}| \pi - \pi_0|^2 = \frac{k_2}{2}(|b|+1)^2|\pi_0|^2 > c
\]
by (\ref{inequality:c}), which contradicts $(R,\Omega) \in V^{-1}([0,c])$. Hence, $b>0$.   If $b>1$, then equation (\ref{abEquations:a}) implies $a^2 <1$, but equation (\ref{abEquations:b}) implies $b^2 < a^2$, implying $b^2 <1$.  This cannot be  compatible with $b>1$. Hence, $b>1$ is ruled out. Similarly,  $0 <b<1$ can be ruled out.  Hence, $b=1$, which implies $\pi = \pi_0$ contradicting (\ref{contra:eqn}). Thus, when $E_ 0= |\pi_0|^2/2I_3$, there are no critical points of $V$ in $V^{-1}([0,c])\backslash V^{-1}(0)$. 

Suppose $E_ 0 \neq |\pi_0|^2/2I_3$. We  analyze equations (\ref{abEquations}) using a continuity argument.  At $a^2=1$,  (\ref{abEquations:a}) implies  $b= 0$ or $1$, neither of which satisfies  (\ref{abEquations:b}) at $a^2=1$. Thus, by continuity  there exists a number $\delta$ with $0 <\delta <1$ such that  for any $a$ with $|a^2 -1|<\delta$ there is no number $b$  satisfying  both  (\ref{abEquations:a}) and (\ref{abEquations:b}). Hence,  $ |a^2-1|\geq\delta$. We now shrink the number $c$ such that it  not only satisfies (\ref{inequality:c}) but also $c< k_0\delta^2 /4$. For such a number $c$, we have
\[
V(R,\Omega) \geq \frac{k_0}{4}\| R^TR-I\|^2 =\frac{k_0}{4}\| (a^2-1)uu^T\|^2  \geq \frac{k_0}{4}\delta^2> c,
\]
which contradicts $(R,\Omega) \in V^{-1}([0,c])$.  Hence, when $E_ 0\neq  |\pi_0|^2/2I_3$, there are no critical points of $V$ in $V^{-1}([0,c])\backslash V^{-1}(0)$ for some $c>0$.

Therefore,  there exists a number $c>0$ such that $V^{-1}(0)$ is the set of all critical points of $V$ in $V^{-1}([0,c])$.
\end{proof}
\end{lemma}

Consider the dynamics
\begin{subequations}\label{rigid:modified}
\begin{align}
\dot R &= R\hat \Omega -k_0R(R^TR - I) - k_2 (\pi(R,\Omega) - \pi_0) \Omega^T \mathbb I, \label{rigid:modified:R} \\
\dot \Omega &= \mathbb I ^{-1} ( (\mathbb I \Omega) \times \Omega ) - k_1 (E(\Omega)-E_0) \mathbb I \Omega - k_2\mathbb IR^T (\pi(R,\Omega) - \pi_0), \label{rigid:modified:Omega}
\end{align}
\end{subequations}
which correspond to (\ref{our:new:dyn:sys}). 
From Theorem \ref{theorem:general} and Lemmas   \ref{lemma:1:rigid}  and \ref{lemma:2:rigid} comes the following theorem: 
\begin{theorem}
There is a number $c>0$ such that every trajectory of (\ref{rigid:modified}) starting from a point in $V^{-1}([0,c])$ stays in $V^{-1}([0,c])$ for all $t\geq 0$ and asymptotically converges to the  set 
\[
V^{-1}(0) = \{ (R,\Omega) \in \mathbb R^{3\times 3} \times \mathbb R^3 \mid R \in \operatorname{SO}(3), \, E(\Omega) = E_0, \, \pi(R,\Omega) = \pi_0\}
\]
 as $t \rightarrow \infty$, where the function $V$ is defined in (\ref{def:V:rigid}). Furthermore, $V^{-1}(0)$ is an invariant set of   both  (\ref{RigidBody}) and (\ref{rigid:modified}).
\end{theorem}

\subsection{The Kepler Problem} 
\label{subsection:two:body}
The  two-body dynamics in the Kepler problem are given in the usual barycentric coordinates by
\begin{subequations}\label{Kepler}
\begin{align}
\dot x &= v,  \label{Kepler:a} \\
\dot v &= -\mu\frac{x}{|x|^3},\label{Kepler:b}
\end{align}
\end{subequations}
where $x \in \mathbb R^3_0:= \mathbb R^3\backslash \{(0,0,0)\}$ is the position vector, $v \in \mathbb R^3$ is the velocity vector and $\mu$ is the gravitational parameter.  Define two  functions $L: \mathbb R^3 \times \mathbb R^3 \rightarrow \mathbb R^3$ and $A : \mathbb R^3_0 \times \mathbb R^3 \rightarrow \mathbb R^3$ by
\begin{align}
L(x,v) &= x \times v, \label{def:L}\\
A(x,v) &= v \times (x \times v) - \mu \frac{x}{|x|},\label{def:A}
\end{align}
where $L$ is called  the angular momentum vector and $A$ is called the Laplace-Runge-Lenz vector. It is known that both $L$ and $A$ are first integrals of the two-body dynamics (\ref{Kepler}) and they are orthogonal to each other, i.e., 
\[
L(x,v) \perp A(x,v)
\]
for all $(x,v) \in \mathbb R^3_0 \times \mathbb R^3$.  The energy function 
\[
E(x,v) = \frac{1}{2}|v|^2 - \frac{\mu}{|x|}
\]
 satisfies 
\begin{equation}\label{Relation}
|A(x,v)|^2 = \mu^2 + 2E(x,v)|L(x,v)|^2
\end{equation}
for all $(x,v) \in \mathbb R^3_0 \times \mathbb R^3$, implying that the energy $E$ is also a first integral of the two-body dynamics (\ref{Kepler}).
 It is also known that a non-degenerate elliptic Keplerian orbit is uniquely determined by a pair $(L,A)$ that satisfies $L \perp A$, $|L| \neq 0$ and $|A| < \mu$ \cite{ChChMa}.

Fix  a non-degenerate elliptic Keplerian orbit, i.e., a  pair of vectors $(L_0, A_0)$ 
that satisfies 
\[
L_0 \perp A_0, \quad |L_0| \neq 0, \quad |A_0| < \mu.
\]
Define a function $V : \mathbb R^3_0 \times  \mathbb R^3 \rightarrow \mathbb R$ by
\begin{equation}\label{def:V:Kepler}
V(x,v) = \frac{k_1}{2}|L(x,v) - L_0|^2 + \frac{k_2}{2} |A(x,v) - A_0|^2
\end{equation}
for $(x,v) \in \mathbb R^3_0 \times \mathbb R^3$, where $k_1>0$ and $k_2>0$. Notice that
\[
V^{-1}(0) = \{ (x,v) \in \mathbb R^3_0 \times \mathbb R^3 \mid L(x,v) = L_0, A(x,v) = A_0\},
\]
which is the non-degenerate Keplerian elliptic orbit whose angular momentum vector and Laplace-Runge-Lenz vector are $L_0$ and $A_0$, respectively. 
\begin{lemma}\label{lemma:gradient:Kepler}
The gradient $(\nabla_x V, \nabla_v V) \in \mathbb R^3 \times \mathbb R^3$ of the function $V$ defined in (\ref{def:V:Kepler}) is given by
\begin{align*}
\nabla_x V&= k_1  v \times \Delta L + k_2 \left (  v \times ( \Delta A \times v ) - \frac{\mu}{|x|} \Delta A + \frac{\mu}{|x|^3} x x^T \Delta A \right ), \\
\nabla_v V  &=k_1 \Delta L \times x + k_2 ( (x \times v) \times \Delta A + x \times ( v \times  \Delta A)),
\end{align*}
where $\Delta L = L(x,v) - L_0$ and $\Delta A = A(x,v) - A_0$.
\end{lemma}

The following lemma shows that the function $V$ defined in (\ref{def:V:Kepler}) satisfies assumptions A1 and A2 stated in \S\ref{section:theory}.
\begin{lemma}\label{ShortLemma:Kepler}
1. The function $V$   satisfies
\[
\langle (\nabla _x V, \nabla_vV), (v, -\mu x /|x|^3) \rangle =0.
\]
2. For any number $c$ satisfying
\begin{equation}\label{range:c:Kepler}
0< c < \min \{k_1 |L_0|^2/2, k_2  (\mu - |A_0|)^2/2 \},
\end{equation}
the set $V^{-1} ([0,c])$ is a compact set in $\mathbb R^3_0 \times \mathbb R^3$.
\begin{proof}
The first fact is a straightforward calculation using the previous Lemma.
For the second, the essential idea is that the fibers of $V$ are homeomorphic to circles, corresponding to the elliptic orbits, and are therefore compact.  For a detailed proof of the second statement,  refer to Corollary 2.2 in \cite{ChChMa}. 
 \end{proof}
 \end{lemma}

The following lemma shows that the function $V$ defined in (\ref{def:V:Kepler}) satisfies assumption A3 stated in \S\ref{section:theory}.

\begin{lemma}\label{LongLemma:Kepler}
For any number $c$ satisfying (\ref{range:c:Kepler})  the set of all critical points of $V$ in $V^{-1}([0,c])$ is equal to $V^{-1}(0)$.
\begin{proof}
Choose an arbitrary number $c$ satisfying (\ref{range:c:Kepler}). 
 Let $(x,v)$ be an arbitrary  critical point of $V$ in $V^{-1}([0,c])$. For notational convenience, let us write
 \[
 L = L(x,v), \quad A = A(x,v)
 \]
 suppressing the dependence on $(x,v)$. By Lemma \ref{lemma:gradient:Kepler},  the critical point $(x,v)$ satisfies 
\begin{subequations}\label{dV:zero}
\begin{align}
0&= k_1  v \times \Delta L + k_2 \left (  v \times ( \Delta A \times v ) - \frac{\mu}{|x|} \Delta A + \frac{\mu}{|x|^3} x x^T \Delta A \right ),  \label{dxV:zero}\\
0&=k_1 \Delta L \times x + k_2 ( (x \times v) \times \Delta A + x \times ( v \times  \Delta A)). \label{dvV:zero}
\end{align}
\end{subequations}
 If $|L| = 0$, then $V(x,v) \geq k_1 |L_0|^2/2 > c$, contradicting $(x,v) \in V^{-1}([0,c])$. Hence, $|L| \neq 0$, which together with (\ref{def:L}) implies that the three vectors $x,v, L$ form a basis for $\mathbb R^3$. The dot product of (\ref{dvV:zero}) with $x$ yields
\[
0 = x \cdot ( (x \times v) \times \Delta A) = \Delta A \cdot (x \times L),
\]
so there are numbers $a$ and $b$ such that 
\begin{equation}\label{DA:xL}
\Delta A = a x + bL. 
\end{equation}
Substitution of (\ref{DA:xL}) into (\ref{dV:zero}) gives
\begin{align*}
0 &= v \times \left ( k_1 \Delta L + k_2 \left ( aL - b v \times L  + \frac{b\mu}{|x|}x\right )\right ),\\
0 &= (k_1 \Delta L  + k_2 (2aL  -b v \times L) ) \times x.
\end{align*}
It follows that there are numbers $d$ and $f$ such that
\begin{subequations}\label{DL:vLx}
\begin{align}
k_1 \Delta L + k_2 \left ( aL - b v \times L  + \frac{b\mu}{|x|}x\right ) =& dv,\label{DL:vLx:1}\\
k_1 \Delta L + k_2 (2aL  -b v \times L)  = &fx. \label{DL:vLx:2}
\end{align}
\end{subequations}
From (\ref{DL:vLx}), we obtain 
\[
 \left ( \frac{bk_2 \mu}{|x|} + f \right )x - dv  - ak_2 L = 0.
\]
By linear independence of $\{x,v, L\}$, 
\[
a = 0, \quad d = 0,\quad f = -bk_2\mu/|x|.
\]
Substitution of these into (\ref{DA:xL}) and (\ref{DL:vLx:2}) gives
\begin{align*}
\Delta A = bL, \quad \Delta L =   \frac{bk_2}{k_1}A,
\end{align*}
where we have used the definition of $A$ given in (\ref{def:A}). Hence,
\begin{equation}\label{A0:L0:linear:combo}
A_0 = A - bL, \quad L_0 = L - \frac{bk_2}{k_1}A.
\end{equation}
From (\ref{A0:L0:linear:combo}) and  the orthogonality  $A_0 \perp L_0$ and $A \perp L$,  it follows that
\[
0 = A_0 \cdot L_0 = -b \left (  |L|^2 + \frac{k_2}{k_1}|A|^2\right ).
\]
Since $|L| \neq 0$, and recalling that $k_{1}>0$ and $k_{2}>0$, we have $b=0$.  Substitution of $b=0$ into (\ref{A0:L0:linear:combo}) yields
\[
L = L_0, \quad A = A_0,
\]
which implies $(x,v) \in V^{-1}(0)$. Thus, every critical point of $V$ in $V^{-1}([0,c])$ is contained in $V^{-1}(0)$. 

Since 0 is the minimum value of $V$, every point in $V^{-1}(0)$ is a critical point of $V$. Therefore, the set of all critical points of $V$ in $V^{-1}([0,c])$ is $V^{-1}(0)$. 
\end{proof}
\end{lemma}
Choose a non-degenerate Keplerian elliptic orbit and let $(x_0, v_0)$ be a point on the orbit. Set
\[
L_0 = L(x_0,v_0), \quad A_0 = A(x_0, v_0)
\]
to be the angular momentum vector and the Laplace-Runge-Lenz vector of the orbit, respectively.  Consider the dynamics:
\begin{subequations}\label{Kepler:modified}
\begin{align}
\dot x &= v -k_1  v \times \Delta L - k_2 \left (  v \times ( \Delta A \times v ) - \frac{\mu}{|x|} \Delta A + \frac{\mu}{|x|^3} x x^T \Delta A \right ),  \label{Kepler:modified:a} \\
\dot v &= -\mu\frac{x}{|x|^3} - k_1 \Delta L \times x - k_2 ( (x \times v) \times \Delta A + x \times ( v \times  \Delta A)),\label{Kepler:modified:b}
\end{align}
\end{subequations}
where $\Delta L = L(x,v) - L_0$ and $\Delta A = A(x,v) - A_0$, which correspond to (\ref{our:new:dyn:sys}). From Theorem \ref{theorem:general} and Lemmas \ref{ShortLemma:Kepler}  and \ref{LongLemma:Kepler} comes the following theorem:

\begin{theorem}
For any $c>0$ satisfying (\ref{range:c:Kepler}),  every trajectory of (\ref{Kepler:modified}) starting from a point in $V^{-1}([0,c])$ stays in $V^{-1}([0,c])$ for all $t\geq 0$ and   asymptotically converges to the  set 
\[
V^{-1}(0) = \{ (x,v) \in \mathbb R^3_0 \times \mathbb R^3 \mid L(x,v) = L_0, A(x,v) = A_0\}
\]
 as $t \rightarrow \infty$, where the function $V$ is defined in (\ref{def:V:Kepler}). Furthermore, $V^{-1}(0)$ is an invariant set of   both  (\ref{Kepler}) and (\ref{Kepler:modified}). 
 \end{theorem}

\subsection{A Perturbed Kepler Problem with Rotational Symmetry}
\label{subsection:pert:Kepler}

Consider a perturbed Kepler problem with rotational symmetry whose equations of motion are given by
\begin{subequations}\label{Pert:Kepler}
\begin{align}
\dot x &= v,  \label{Pert:Kepler:a} \\
\dot v &= -U^\prime(|x|) \frac{x}{|x|},\label{Pert:Kepler:b}
\end{align}
\end{subequations}
where $x \in \mathbb R^3_0:= \mathbb R^3\backslash \{(0,0,0)\}$ is the position vector, $v \in \mathbb R^3$ is the velocity vector, and $U : (0, \infty) \rightarrow \mathbb R$ is the potential function  that depends only on the radial distance from the origin. The total energy $E: \mathbb R^3_0 \times \mathbb R^3 \rightarrow \mathbb R$ and the angular momentum vector $L : \mathbb R^3_0 \times \mathbb R^3 \rightarrow \mathbb R^3$ are defined by
\begin{align}
E(x,v) &= \frac{1}{2}|v|^2 + U(|x|), \label{Pert:Kep:energy}\\
L(x,v) &= x \times v \label{Pert:Kep:ang:mom}
\end{align}
and they are conserved quantities of the  dynamics (\ref{Pert:Kepler}). 
Take any point $(x_0, v_0) \in \mathbb R^3_0 \times \mathbb R^3$ such that
\[
x_0 \times v_0 \neq 0.
\]
Let
\[
E_0 = E(x_0, v_0), \quad L_0 = L(x_0, v_0) \neq 0.
\]
Define a function $V: \mathbb R^3_0 \times \mathbb R^3 \rightarrow \mathbb R$ by
\[
V(x,v) = \frac{k_1}{2}  |E(x,v)-E_0|^2 + \frac{k_2}{2} |L(x,v)-L_0|^2
\]
with $k_1>0$ and $k_2>0$. Then,
\[
V^{-1}(0) = \{ (x,v) \in \mathbb R^3_0 \times \mathbb R^3 \mid E(x,v) = E_0, L(x,v) = L_0\}.
\]
The gradient $(\nabla_x V, \nabla_v V)$ of $V$ is given by
\begin{align*}
\nabla_x V &= k_1 \Delta E U^\prime(|x|) \frac{x}{|x|} + k_2 v\times \Delta L,\\
\nabla_v V &= k_1 \Delta E v  + k_2  \Delta L \times x,
\end{align*}
where $\Delta E = E(x,v) - E_0$ and $\Delta L = L(x,v) - L_0$. Trivially, $V$ satisfies (\ref{Vdot:X}), i.e.
\begin{equation}\label{pert:VX0}
\langle (\nabla_xV, \nabla_vV), (v, -U^\prime(|x|)x/|x|\rangle = 0
\end{equation}
for all $(x,v) \in \mathbb R_0^3 \times \mathbb R^3$.
The modified dynamics, which correspond to (\ref{our:new:dyn:sys}), are computed as 
\begin{subequations}\label{Pert:Kepler:Mod}
\begin{align}
\dot x &= v - k_1 \Delta E U^\prime(|x|) \frac{x}{|x|} - k_2 v\times \Delta L,  \label{Pert:Kepler:Mod:a} \\
\dot v &= -U^\prime(|x|) \frac{x}{|x|} -k_1 \Delta E v  - k_2  \Delta L \times x.\label{Pert:Kepler:Mod:b}
\end{align}
\end{subequations}
\begin{theorem}\label{theorem:main:pert:Kep}
Suppose that $V^{-1}(0)$ is compact and there is no common solution $r>0$ to the following two equations:
\begin{align}
E_0 &= \frac{1}{2}rU^\prime(r) + U(r), \label{E0:cond:Pert}\\
|L_0|^2 &= r^3 U^\prime (r).\label{L0:cond:Pert}
\end{align}
Then,  assumptions A2 and A3 hold and  there is a number $c>0$ such that every trajectory of (\ref{Pert:Kepler:Mod}) starting in $V^{-1}([0,c])$  remains in $V^{-1}([0,c])$ for all $t\geq 0$ and asymptotically converges to $V^{-1}(0)$ as $t\rightarrow \infty$. 
\begin{proof}
Define a function  $f: \mathbb R^3_0 \times \mathbb R^3 \rightarrow \mathbb R \times \mathbb R^3$  by
\[
f(x,v) = \begin{bmatrix}
E(x,v) \\ L(x,v)
\end{bmatrix}.
\]
Then,
\[
Df(x,v)^T = \begin{bmatrix}
U^\prime(|x|) \frac{x}{|x|} & \hat v  \\
 v & -\hat x
\end{bmatrix},
\]
where the over-hat symbol $\wedge$ denotes the hat map defined in (\ref{def:hat:map}). We want to show that the $6\times 4$ matrix $Df(x,v)^T$ is one-to-one for all $(x,v) \in V^{-1}(0)$.  Fix an arbitrary point $(x,v) \in V^{-1}(0)$. It follows
\begin{align}
E_0 &= \frac{1}{2}|v|^2 + U(|x|),\label{thm:E0}\\
L_0 &= x \times v \neq 0. \label{thm:L0}
\end{align}
Take any point $(a,w) \in \mathbb R \times \mathbb R^3$ from the kernel of $Df(x,v)^T$. Then, 
\begin{subequations}\label{DfT:ker}
\begin{align}
0& = aU^\prime(|x|)\frac{x}{|x|} + v \times w, \label{DfT:ker:a}\\
0&=a  v - x \times w. \label{DfT:ker:b}
\end{align}
\end{subequations}

Suppose $a\neq 0$. Taking the inner product of (\ref{DfT:ker:a}) with $x$ and of (\ref{DfT:ker:b}) with $v$, we obtain 
\begin{align*}
0&= a U^\prime(|x|)  |x| + L_0 \cdot w,\\
0&= a |v|^2 + L_0 \cdot w,
\end{align*}
from which it follows that
\begin{equation}\label{Uprime:v2}
|x|  U^\prime(|x|) = |v|^2.
\end{equation}
Taking the inner product of (\ref{DfT:ker:b}) with $x$, we get  $x\cdot v=0$ which implies
\begin{equation}\label{L0:norm}
|L_0| = |x|\cdot |v|.
\end{equation}
From (\ref{thm:E0}), (\ref{Uprime:v2}) and (\ref{L0:norm}), we obtain
\begin{align}
E_0 &= \frac{1}{2}|x| U^\prime(|x|) + U(|x|), \label{E0:eqn:Up}\\
|L_0|^2 &= |x|^3 U^\prime(|x|).\label{L0:norm:eqn}
\end{align}
By hypothesis, there cannot be any $x \in \mathbb R_0^3$ that satisfies both (\ref{E0:eqn:Up}) and (\ref{L0:norm:eqn}). Hence, we cannot have $a\neq 0$. 

Substitute $a=0$ into (\ref{DfT:ker}). It follows that $w$ is parallel to $x\times v$. Hence, there is a number $b$ such that $w = b L_0$. Substituting this in (\ref{DfT:ker:b}) yields $ b x \times L_0 = 0$. Taking the cross product of this with  $x$ yields  $b|x|^2L_0 = 0$ since $x\cdot L_0=0$.  Since $x\neq 0$ and $L_0 \neq 0$, we have $b=0$, so $w = 0$. It follows that $(a,w) = (0,0)$, which implies that  $Df(x,v)^T$ is one-to-one for all $(x,v) \in V^{-1}(0)$. In other words, $Df(x,v)$ is onto for all $(x,v) \in V^{-1}(0)$. Hence, the conclusion of the theorem follows from Lemma \ref{kuhl:lemma}, equation (\ref{pert:VX0}), and Theorem \ref{kuhl:theorem}.

\end{proof}
\end{theorem}

\begin{remark}
Consider a special case in which the potential function $U(r)$ is of the form
\begin{equation}\label{pert:kep:pot}
U(r) = -\frac{\mu}{r} - \frac{\delta}{r^3},
\end{equation}
where $\mu>0$ and $\delta>0$.
Then equations (\ref{E0:cond:Pert}) and (\ref{L0:cond:Pert}) become
\begin{align}
E_0 &= -\frac{\mu}{2r} + \frac{\delta }{2r^3}, \label{E0:Birk} \\
|L_0|^2 &= \mu r + \frac{3\delta }{r}. \label{L0:Birk}
\end{align}
Given $E_0$ and $L_0$, it is then easy to check if there is no common solution $r>0$ to  (\ref{E0:Birk}) and (\ref{L0:Birk}). 
\end{remark}

\section{Simulations}

\subsection{The Free Rigid Body}

Consider the free rigid body dynamics in \S\ref{subsection:rigid:body} with the moment of inertia matrix $\mathbb I = \operatorname{diag} (3,2,1)$ and  the initial condition
\begin{equation}\label{IC:rigid}
R(0) = I, \quad \Omega (0) = (1,1,1).
\end{equation}
The values of the energy $E$ and the spatial angular momentum vector $\pi = (\pi_1, \pi_2, \pi_3)$ corresponding to the initial condition are
\[
E (0)= 3, \quad \pi(0)  = (3,2,1).
\]
The period $T_\Omega$ of the trajectory of the body angular velocity vector $\Omega(t) $ is computed approximately to be $T_\Omega = 6.4227$.  

We integrate the dynamics over the time interval $[0, 10^3] = [0, 155.7T_\Omega]$ with  step size $\Delta t = 10^{-4}$, using the  following four integration methods: a feedback integrator with the Euler scheme, a projection method with the Euler scheme,  a splitting method with three rotations splitting, and the ordinary Euler method.  The feedback integrator with the Euler scheme denotes  the Euler method   applied to the modified free rigid dynamics (\ref{rigid:modified}) with the following values of the parameters $k_0$, $k_1$, and $k_2$ 
\[
k_0 = 50, \quad k_1 = 100, \quad k_2 = 50.
\]
The projection method is the standard one explained on pp.110--111 in \cite{Hair}. In order to solve  constraint equations for projection at each step of integration in the projection method, we use  the Matlab command \textit{fsolve}  with the parameter \textit{TolFun}, which is termination tolerance on the function value,  set equal to  $10^{-4}$, which is the same as the integration step size $\Delta t$. The splitting method is the one explained on pp.284--285 in \cite{Hair}. The three of the projection method, the splitting method and the ordinary Euler method are applied to the original free rigid body dynamics (\ref{RigidBody}). 

\begin{figure}[t]
 \includegraphics[scale = 0.65]{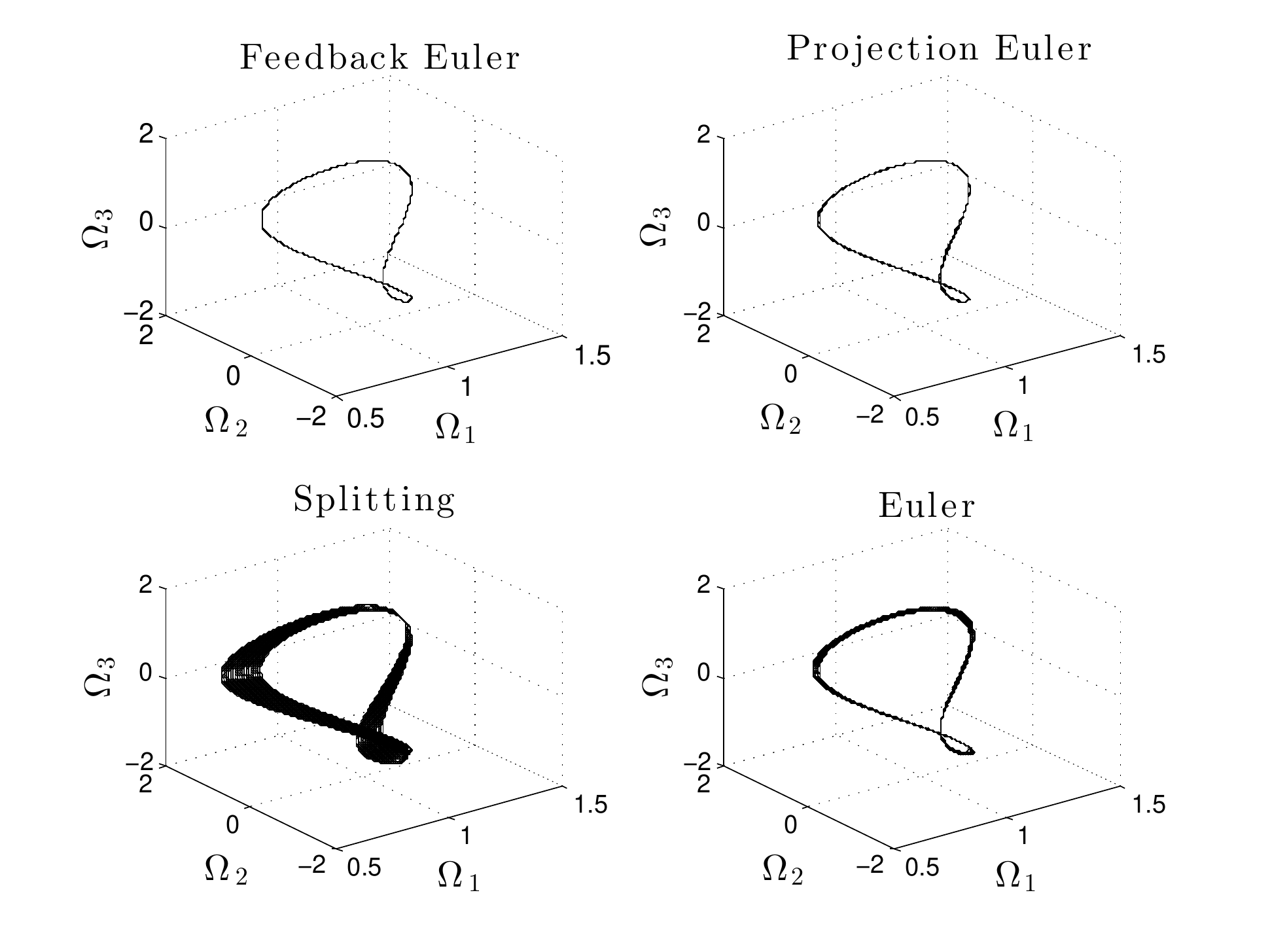}
\caption{The trajectories of the body angular velocity $\Omega (t) = (\Omega_1(t), \Omega_2(t), \Omega_3(t))$, $0 \leq t \leq 1000$, of the free rigid body dynamics generated by  four different methods with step size $\Delta t = 10^{-4}$: a feedback integrator with the Euler scheme, the standard projection method with the Euler scheme, a three rotations splitting method and the usual Euler method.}
\label{figure:Rigid:Omega}
\end{figure}

\begin{figure}[h]
 \includegraphics[scale = 0.5]{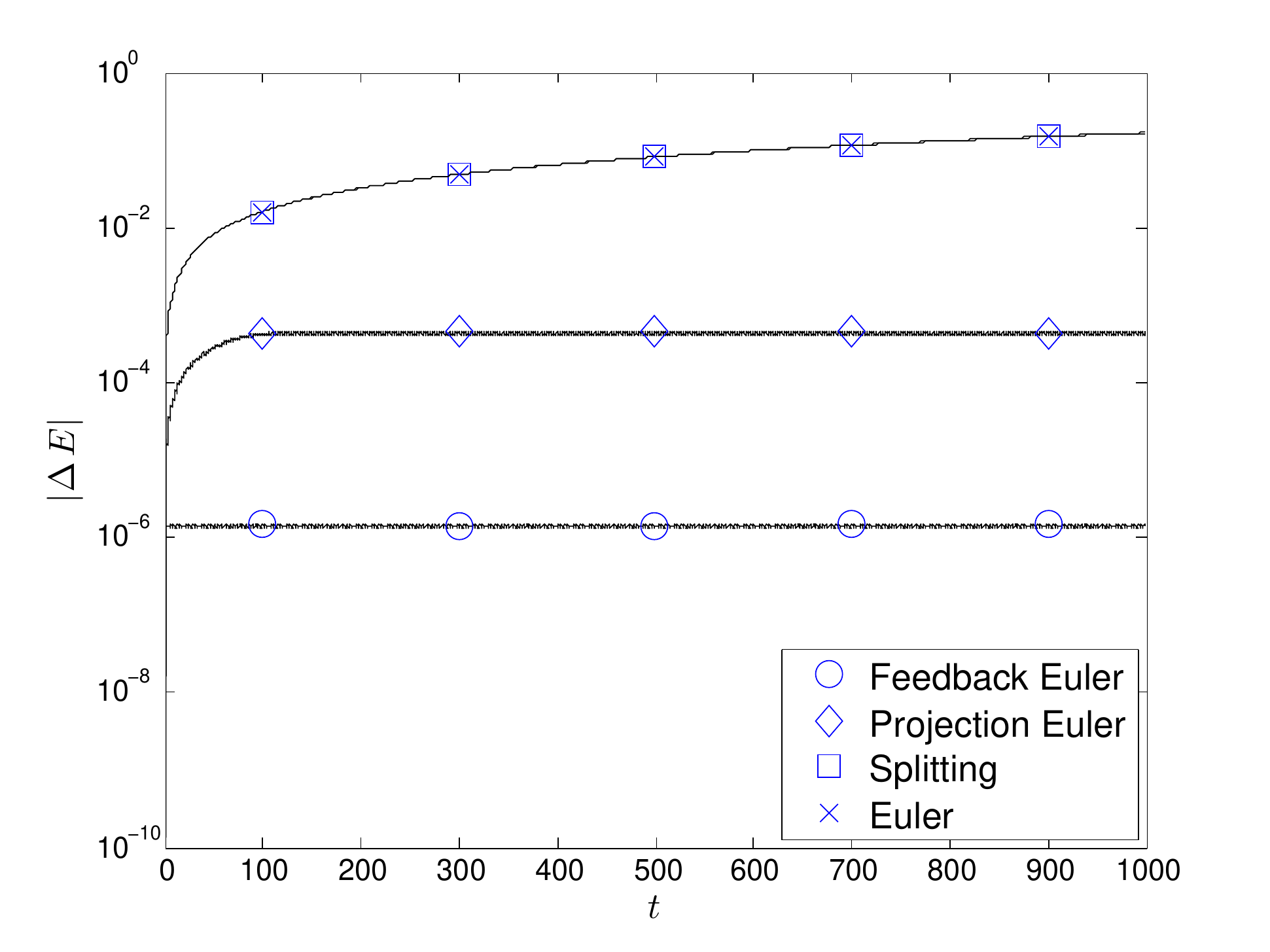}
\caption{The trajectories of the energy error $|\Delta E(t) | = |E(t) - E(0)|$, $0 \leq t \leq 1000$, of the free rigid body dynamics  generated by  four different methods with step size $\Delta t = 10^{-4}$: a feedback integrator with the Euler scheme ($\circ$), the standard projection method with the Euler scheme ($\diamond$), a three rotations splitting method ($\square$) and the usual Euler method ($\times$).}
\label{figure:Rigid:Energy}
\end{figure}

\begin{figure}[h]
 \includegraphics[scale = 0.5]{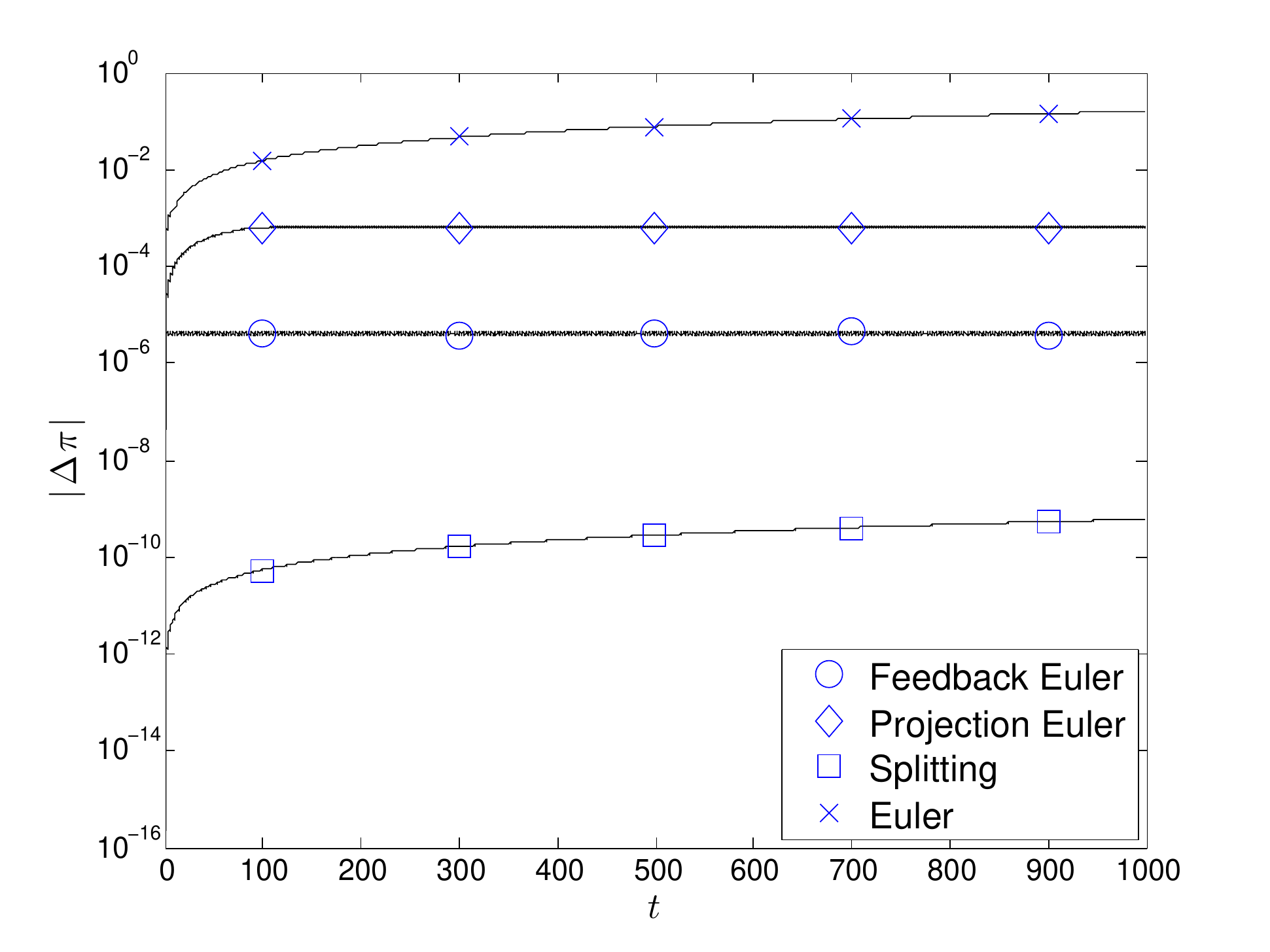}
\caption{The trajectories of the spatial angular momentum error  $|\Delta \pi (t) | = |\pi (t) - \pi (0)|$, $0 \leq t \leq 1000$, of the free rigid body dynamics  generated by  four different methods with step size $\Delta t = 10^{-4}$: a feedback integrator with the Euler scheme ($\circ$), the standard projection method with the Euler scheme ($\diamond$), a three rotations splitting method ($\square$) and the usual Euler method ($\times$).}
\label{figure:Rigid:Momentum}
\end{figure}

\begin{figure}[h]
 \includegraphics[scale = 0.5]{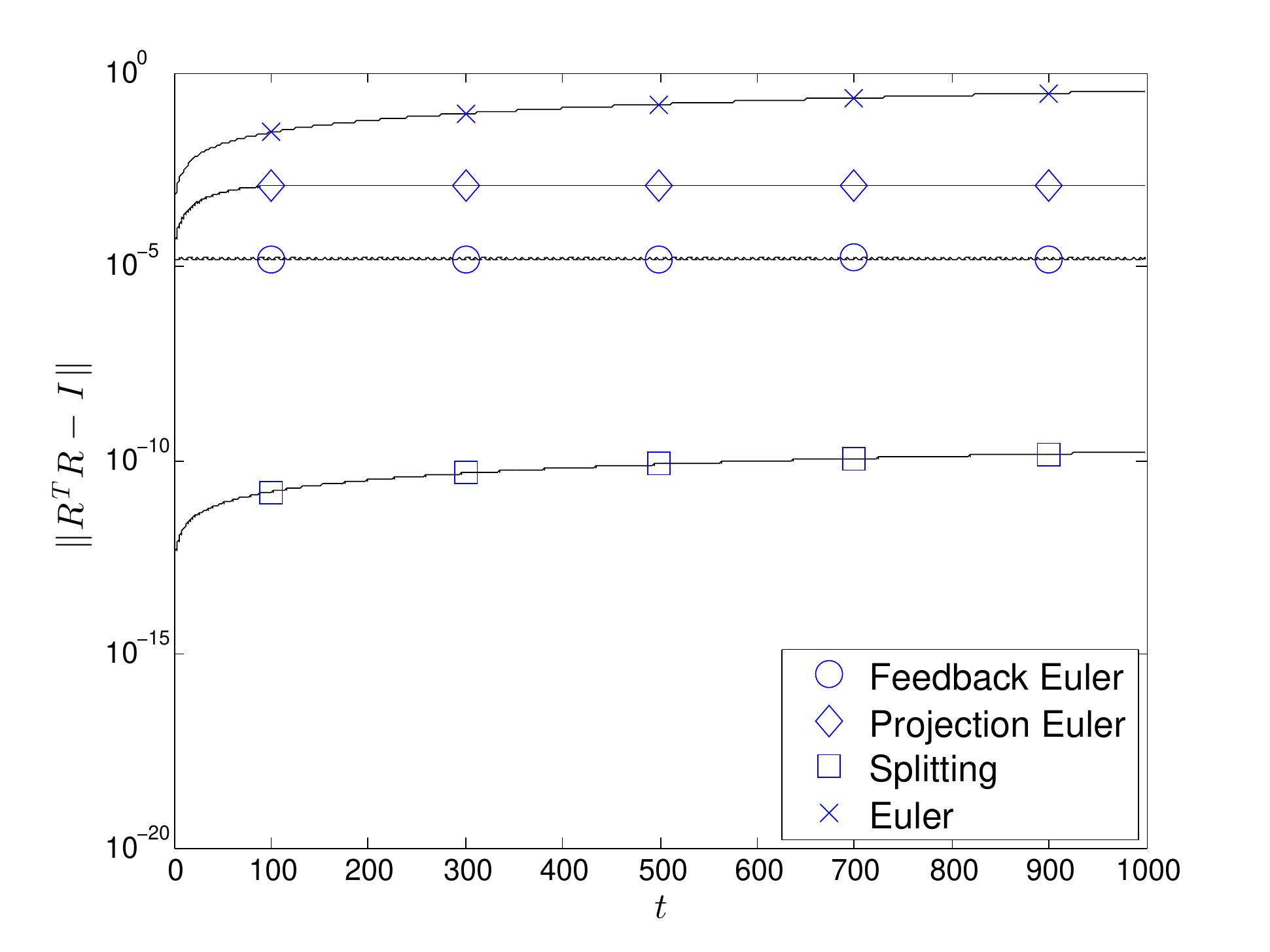}
\caption{The trajectories of the deviation   $\| R (t)^TR(t) - I \| $ of the rotation matrix $R(t)$ from ${\operatorname{SO}(3)}$, $0 \leq t \leq 1000$, of the free rigid body dynamics  generated by  four different methods with step size $\Delta t = 10^{-4}$: a feedback integrator with the Euler scheme ($\circ$), the standard projection method with the Euler scheme ($\diamond$), a three rotations splitting method ($\square$) and the usual Euler method ($\times$).}
\label{figure:Rigid:SO3}
\end{figure}

The trajectories of the body angular velocity vector $\Omega(t)$, the energy error  $|\Delta E (t)| = |E(t) - E(0)|$, the error $|\Delta \pi (t) | = |\pi(t) - \pi(0)|$ in spatial angular momentum, and the deviation $\| R(t)^T R(t) - I\|$ of the rotation matrix $R(t)$ from ${\operatorname{SO}(3)}$ are plotted in Figures \ref{figure:Rigid:Omega}, \ref{figure:Rigid:Energy}, \ref{figure:Rigid:Momentum} and \ref{figure:Rigid:SO3}, respectively. In Figure \ref{figure:Rigid:Omega},  it is observed that the trajectories of $\Omega(t)$ generated by the feedback integrator and the projection method maintain a periodic shape well whereas those by the splitting method and the Euler method drift away significantly from the periodic shape.   In Figure  \ref{figure:Rigid:Energy}, it is observed that the feedback integrator and the projection method keep the energy error sufficiently small whereas the energy errors by the other two methods  increase in time. Although the two  trajectories of  energy error by the splitting method and the  Euler method seem to coincide in Figure \ref{figure:Rigid:Energy}, an examination of the numerical data shows that the energy of the Euler method gets larger than that of the splitting method in time. For example, at $t=1000$, the energy of the Euler method is bigger than that of the splitting method by $1.767\times 10^{-3}$. In Figures \ref{figure:Rigid:Momentum} and \ref{figure:Rigid:SO3}, it is observed  that the feedback method  preserves the spatial angular momentum vector and the  manifold ${\operatorname{SO}(3)}$ sufficiently well.  In terms of computation time, the projection method takes much more time than the others, which is due to the steps of solving the constraint equations for projection. The splitting method is symplectic and of order 2 whereas the other methods are of order 1.  All of these observations  lead us to the conclusion that  the feedback integrator overall has produced the best outcome  in the simulation of the free rigid body dynamics.

\subsection{The Kepler Problem}
\label{subsection:simulation:Kepler}

Consider the Kepler problem in \S\ref{subsection:two:body} with $\mu = 1$ and the initial condition
\[
x(0) = (1,0,0), \quad v(0) = (0,\sqrt{1.8},0).
\]
The corresponding initial values of the angular momentum vector and the Laplace-Runge-Lenz vector  are 
\[
L(0) = (0,0,\sqrt{1.8}), \quad A(0) = (0.8, 0,0).
\]
The  period $T$  and the eccentricity $e$ of the Kepler orbit containing the initial point are
\[
T = 70.2481,\quad e = 0.8.
\]
We integrate the Kepler dynamics over the time interval $[0, 1000T]$ with  step size $\Delta t = 0.005$, using  the following four integration methods: a feedback integrator with the Euler scheme, the standard projection method with the Euler scheme, and two St\"ormer-Verlet schemes.  The feedback integrator with the Euler scheme denotes  the Euler method   applied to (\ref{Kepler:modified})  with $k_1 = 4$ and $k_2 = 2$. The standard projection method is explained  on pp.110--111 in \cite{Hair}.  To solve the constraint equations for projection, we use the Matlab command \textit{fsolve}   with  the parameter \textit{TolFun} set equal to $0.005$, which is the same as  the integration step size $\Delta t$. The two St\"ormer-Verlet schemes are those in (3.4) and (3.5)  on pp. 189--190 in \cite{Hair}, and we call them St\"ormer-Verlet-A and St\"ormer-Verlet-B, respectively,  for convenience. The  St\"ormer-Verlet schemes are symplectic methods of order 2. 

\begin{figure}[h]
 \includegraphics[scale = 0.55]{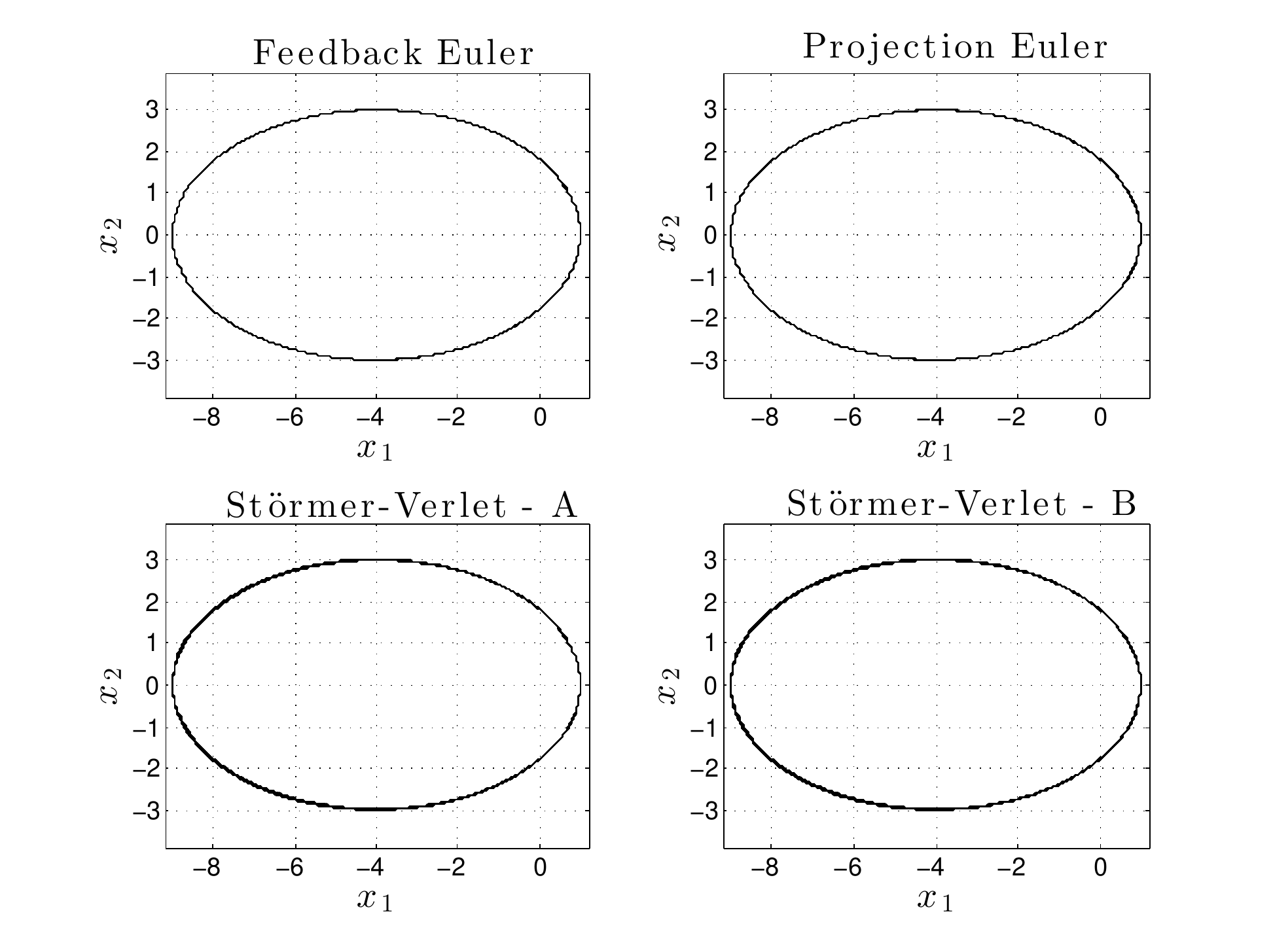}
\caption{The trajectories of the planar orbit $x(t) = (x_1(t), x_2(t),0)$, $0 \leq t \leq 70,248$, in the Kepler problem  generated by four different methods with step size $\Delta t = 0.005$: a feedback integrator with the Euler scheme, the standard projection method with the Euler scheme, and two St\"ormer-Verlet schemes. }
\label{figure:KepEllipse}
\end{figure}

\begin{figure}[h]
 \includegraphics[scale = 0.6]{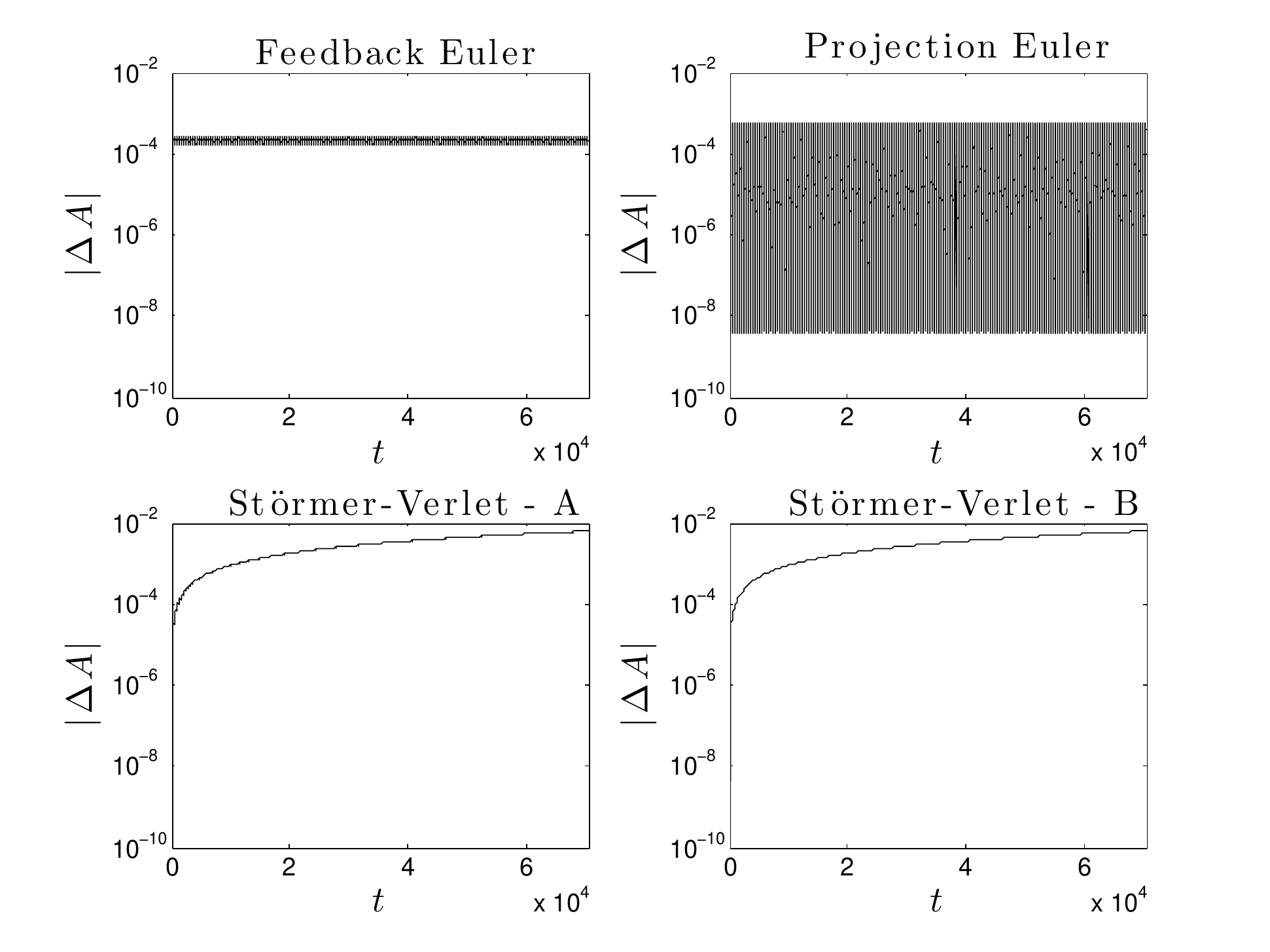}
\caption{The trajectories of the  error $|\Delta A (t)| = |A(t) - A(0)|$, $0 \leq t \leq 70,248$, of the Laplace-Runge-Lenz vector in the Kepler problem  generated by four different methods with step size $\Delta t = 0.005$: a feedback integrator with the Euler scheme, the standard projection method with the Euler scheme, and two St\"ormer-Verlet schemes.}
\label{figure:KepA}
\end{figure}

\begin{figure}[h]
 \includegraphics[scale = 0.55]{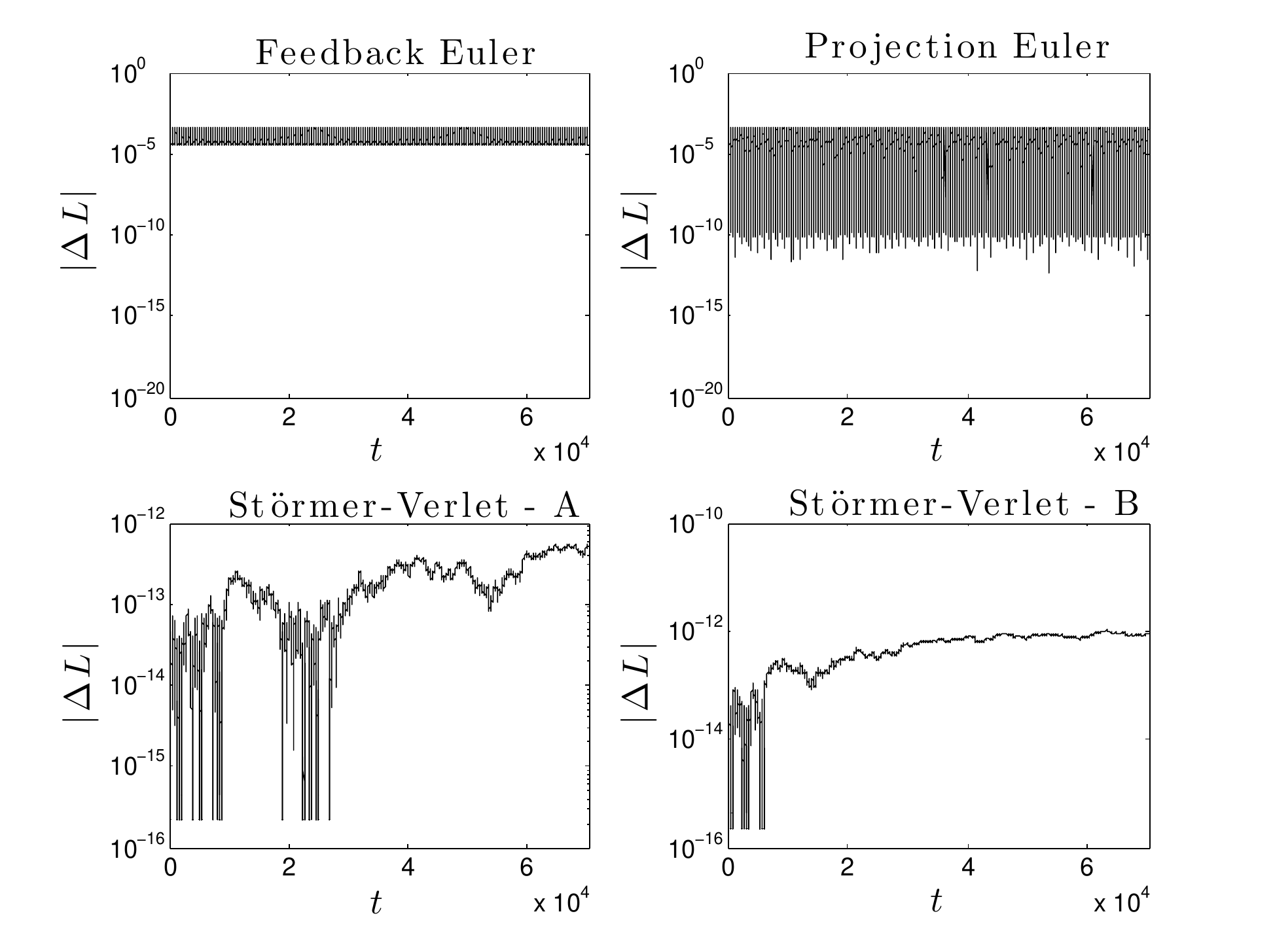}
\caption{The trajectories of the angular momentum  error $|\Delta L (t)| = |L(t) - L(0)|$, $0 \leq t \leq 70,248$, in the Kepler problem generated  by four different methods with step size $\Delta t = 0.005$: a feedback integrator with the Euler scheme, the standard projection method with the Euler scheme, and two St\"ormer-Verlet schemes.}
\label{figure:KepL}
\end{figure}

The trajectories of the  planar orbit $x(t) = (x_1(t), x_2(t),0)$, the error of the Laplace-Runge-Lenz vector, $|\Delta A (t) | = |A(t) - A(0)|$,  and the error of the angular momentum vector,  $|\Delta L(t)| = |L(t) - L(0)|$, are plotted in Figures \ref{figure:KepEllipse}, \ref{figure:KepA} and \ref{figure:KepL}. In Figure \ref{figure:KepEllipse} it is observed that the planar trajectories $x(t) = (x_1(t), x_2(t),0)$ generated by the feedback integrator and the projection method maintain the elliptic shape well whereas those by the St\"ormer-Verlet schemes precess. This can be also verified   in Figure \ref{figure:KepA}, where the feedback integrator and the projection method preserve  the Laplace-Runge-Lenz vector well, but the St\"ormer-Verlet schemes cause the Laplace-Runge-Lenz vector to noticeably precess. In Figure \ref{figure:KepL}, it is observed that the  St\"ormer-Verlet schemes preserve the angular momentum vector exceptionally well in comparison with the other two methods. In Figures \ref{figure:KepA} and \ref{figure:KepL}, we can see that the precision of the feedback integrator is comparable with that of the projection method. However, the feedback integrator takes much less computation time than the projection method.  The feedback integrator and the projection method used here are of order 1, whereas the St\"ormer-Verlet schemes are of order 2.  All of these observations lead us to conclude that the feedback integrator has produced the best result overall.

\subsection{A Perturbed Kepler Problem with Rotational Symmetry} 
Consider the perturbed Kepler problem in \S\ref{subsection:pert:Kepler} with the potential function $U$ given in (\ref{pert:kep:pot}) with $\mu = 1$ and $\delta = 0.0025$, which is the one used in Example 4.3 on p. 111 in \cite{Hair}.  We use the  initial conditions
\[
x(0) = (1-e,0,0), \quad v(0) = (0, \sqrt{(1+e)/(1-e)},0)
\]
with eccentricity $e= 0.6$ as  in \cite{Hair}.  The corresponding values of the energy and the angular momentum vector are
\[
E(0) = -0.5390625, \quad  L (0) = (0,0,0.8).
\]
 We integrate the perturbed Kepler dynamics over the time interval $[0, 200]$ with  step size $\Delta t = 0.03$, just as on  p. 111 in \cite{Hair}, using the  following four integration methods: a feedback integrator with the Euler scheme, the standard projection method  with the Euler scheme, the St\"ormer-Verlet scheme  in (3.4) on p. 189 in \cite{Hair}, and the Matlab command, \textit{ode45}. 
The feedback integrator with the Euler scheme denotes  the Euler method   applied to (\ref{Pert:Kepler:Mod}) with $k_1 = 2$ and $k_2 = 3$, and it is straightforward to verify that the hypotheses in Theorem \ref{theorem:main:pert:Kep} hold true. The other three methods are applied to (\ref{Pert:Kepler}). The Matlab command \textit{fsolve} is used in the projection method with the parameter \textit{TolFun} set equal to $10^{-8}$.  The options of \textit{RelTol} = \textit{AbsTol} = $10^{-10}$ are used for the Matlab integrator, \textit{ode45}, so the result generated by \textit{ode45} can be used as a reference. 

\begin{figure}[h]
 \includegraphics[scale = 0.55]{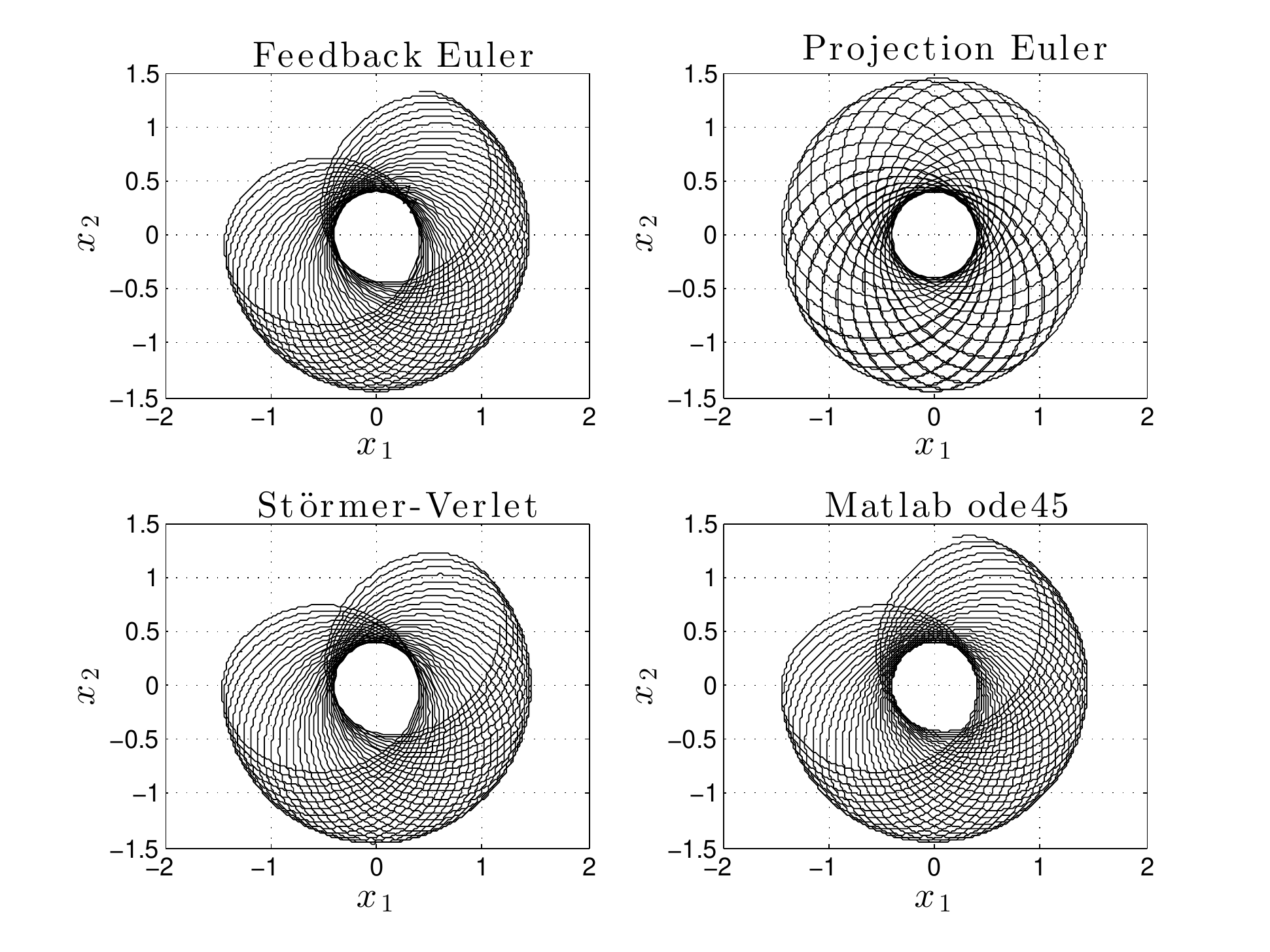}
\caption{The trajectories of the planar orbit $x(t) = (x_1(t), x_2(t),0)$, $0 \leq t \leq 200$, in the perturbed Kepler problem  generated by four different methods: a feedback integrator with the Euler scheme, the standard projection method with the Euler scheme, a St\"ormer-Verlet scheme and the Matlab command \textit{ode45}, where the step size $\Delta t = 0.03$ is used for the first three methods.}
\label{figure:PertKepEllipse}
\end{figure}

\begin{figure}[h]
 \includegraphics[scale = 0.55]{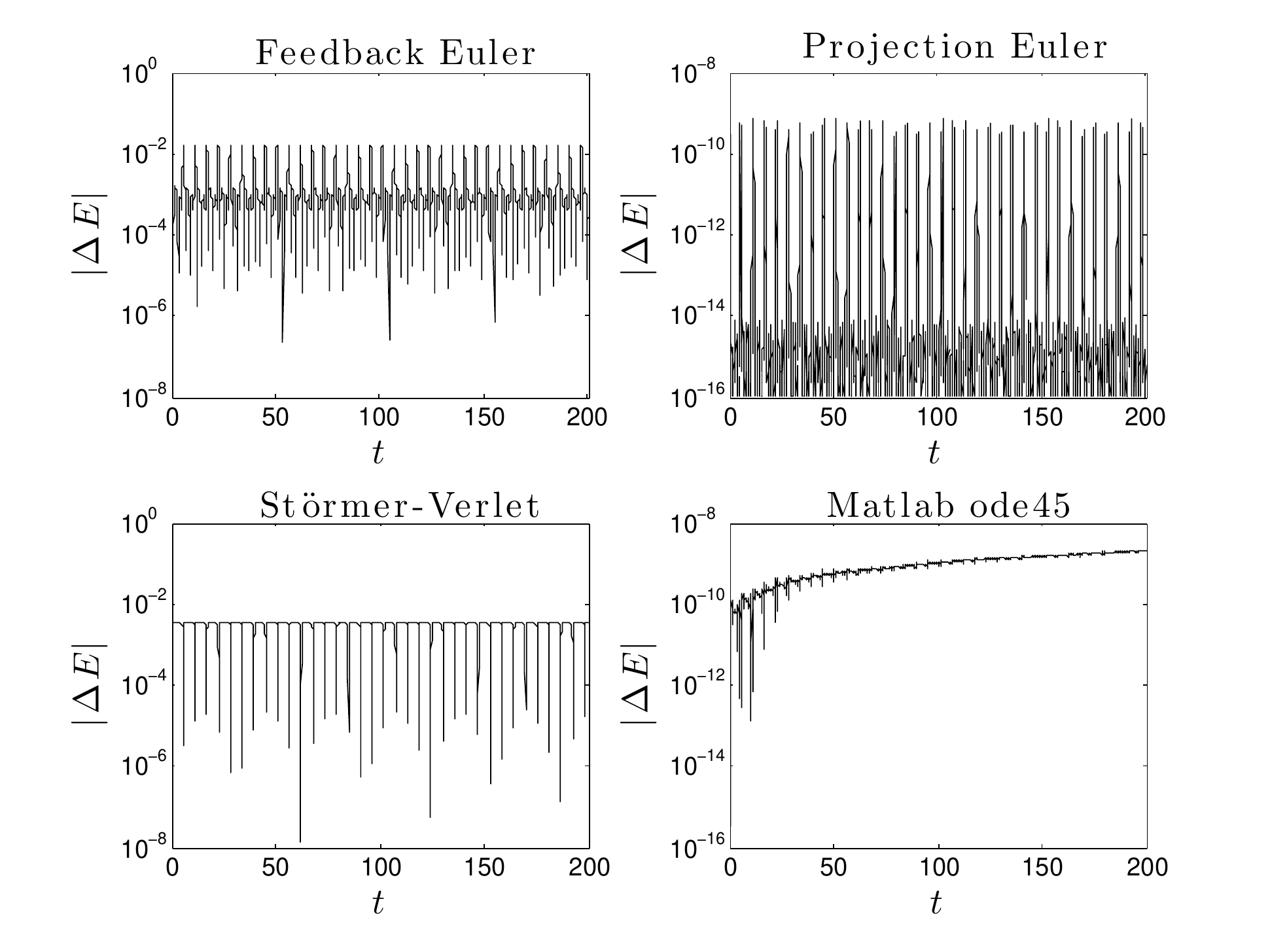}
\caption{The trajectories of the energy error $|E (t)| = |E(t) - E(0)|$, $0 \leq t \leq 200$, in the perturbed Kepler problem  generated by four different methods: a feedback integrator with the Euler scheme, the standard projection method with the Euler scheme, a St\"ormer-Verlet scheme and the Matlab command \textit{ode45}, where the step size $\Delta t = 0.03$ is used for the first three methods.}
\label{figure:PertKepEnergy}
\end{figure}

\begin{figure}[h]
 \includegraphics[scale = 0.55]{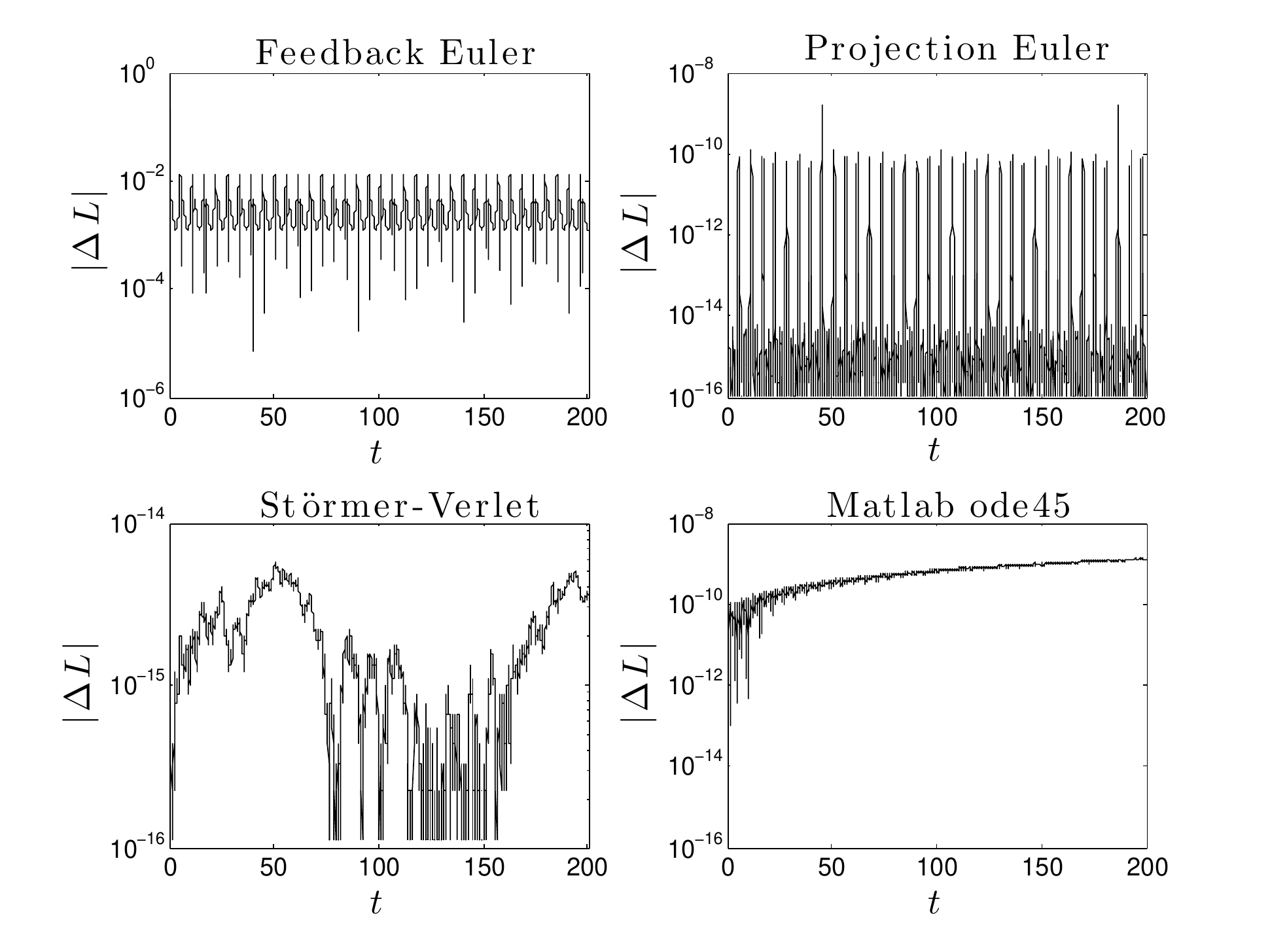}
\caption{The trajectories of the angular momentum error $|\Delta L (t)| = |L(t) - L(0)|$, $0 \leq t \leq 200$, in the perturbed Kepler problem  generated by four different methods: a feedback integrator with the Euler scheme, the standard projection method with the Euler scheme, and a St\"ormer-Verlet scheme and the Matlab command \textit{ode45}, where the step size $\Delta t = 0.03$ is used for the first three methods.}
\label{figure:PertKepL}
\end{figure}

The trajectories of the  planar orbit $x(t) = (x_1(t), x_2(t),0)$, the energy error $|\Delta E (t) | = |E(t) - E(0)|$ and the error $| \Delta L(t)| = |L(t) - L(0)|$ in angular momentum are plotted in Figures \ref{figure:PertKepEllipse}, \ref{figure:PertKepEnergy} and \ref{figure:PertKepL}. 
In Figure \ref{figure:PertKepEllipse} it is observed that the orbits generated by the feedback integrator and the St\"ormer-Verlet scheme are similar to that by \textit{ode45}, but the orbit by the projection method precesses too much which is a very poor result.  The projection method excels only at preserving the energy and the angular momentum as expected in view of the nature of the projection method and the small tolerance parameter value, \textit{TolFun} = $10^{-8}$,  used for the Matlab command, \textit{fsolve}.  In Figure  \ref{figure:PertKepEnergy}, it is observed that the feedback integrator is comparable with the  St\"ormer-Verlet scheme in energy conservation. The feedback integrator also preserves the angular momentum  well in view of the step size $\Delta t = 0.03$, as can be seen in Figure \ref{figure:PertKepL}. The feedback integrator and the projection method used here are of order 1 whereas the St\"ormer-Verlet scheme is of order 2.  From all of these observations,  we conclude that  the feedback integrator  has  produced the best result overall.

\section{Conclusions and Future Work}

We have developed a theory to produce numerical trajectories of a dynamical system on a manifold that  stably remain on the manifold and preserve first integrals of the system. Our theory is not a numerical  integration scheme but rather a modification of   the original dynamics by feedback. The actual numerical integration in our framework can be done  with any usual integrator such as Euler and Runge-Kutta. Our method is successfully applied to the free rigid body,  the Kepler problem and a perturbed Kepler problem with rotational symmetry,  and its excellent performance is demonstrated by simulation studies in comparison with the standard projection method,  two St\"ormer-Verlet schemes and a splitting method via three rotations splitting.

As future work, we plan to   apply our theory to various mechanical  systems with symmetry and non-holonomic systems. We also plan to  carry out a quantitative study of the effect of the parameters in the Lyapunov function on the performance of our method.

\section*{Appendix}

We  show, using results in \cite{KlLo86}, that  any discrete-time dynamical system derived from a one-step numerical integration scheme with uniform step size $h$ for the modified dynamical system (\ref{our:new:dyn:sys}) has an attractor $\Lambda_h$ that contains $V^{-1}(0)$ in its interior and converges to $V^{-1}(0)$ as $h\rightarrow 0+$.  Let us first review some definitions from  \cite{KlLo86}.  Let $A$ and $B$ be nonempty,  compact subsets of $\mathbb R^n$ and $x$ a point in $\mathbb R^n$.  
The  distance between $x$ and $A$ is defined by
\[
\operatorname{dist} (x,A) = \inf\{ |x- a|,  a \in A\}.
\]
The Hausdorff separation of $A$ from $B$ is defined by
\[
H^*(A,B) = \max\{\operatorname{dist}(a,B),  a \in A\}.
\]
The Hausdorff distance between $A$ and $B$ is defined by
\[
H(A,B) = \max \{ H^*(A,B), H^*(B,A)\}.
\]
The Hausdorff distance is a metric on the space of nonempty compact subsets of $\mathbb R^n$.  For $r>0$, let 
\[
S(A,r) = \{ x \in \mathbb R^n \mid \operatorname{dist}(x,A)<r\}
\]
denote an $r$-neighborhood of $A$. 

We say that a nonempty, compact subset $\Lambda$ of $\mathbb R^n$ is uniformly stable for an autonomous dynamical system  if for each $\epsilon >0$ there exists a $\delta = \delta(\epsilon) >0$ such that 
\[
 [x_0 \in S(\Lambda, \delta)  \textup{ and }  t \geq 0 ] \Rightarrow x(t;x_0) \in S(\Lambda, \epsilon), 
\]
where $x(t;x_0)$ is the solution of the given dynamical system with initial condition $x(0) = x_0$.   A set $\Lambda$ is said to be positively invariant for an autonomous  dynamical system if $x(t;x_0) \in \Lambda$ for all $x_0 \in \Lambda$ and $t\geq0$. A nonempty, compact subset $\Lambda$ of $\mathbb R^n$ is called  uniformly asymptotically stable for an autonomous dynamical system if it is  positively invariant  and uniformly stable for the dynamical system, and additionally satisfies the following property:   there is a $\delta_0>0$ and for each $\epsilon>0$ a time $T(\epsilon)>0$ such that 
\[
 [x_0 \in S(\Lambda, \delta_0)  \textup{ and }  t \geq T(\epsilon) ] \Rightarrow x(t;x_0) \in S(\Lambda, \epsilon).
\]
\begin{lemma}\label{lemma:UAS:Lambda}
Suppose that assumptions A1, A2 and  A3 (or  A3$^{\prime}$ instead of A3) stated in \S\ref{section:theory} hold true. Then, the set $V^{-1}(0)$ is uniformly asymptotically stable for the modified dynamical system (\ref{our:new:dyn:sys}).
\begin{proof}
Since the three assumptions are satisfied, the conclusions of Theorem \ref{theorem:general} (or, \ref{theorem:general:more}) hold true.  For convenience, let $\Lambda = V^{-1}(0)$, which is invariant under (\ref{our:new:dyn:sys}) by Theorem \ref{theorem:general} (or, \ref{theorem:general:more}).  Let $c>0$ be the number $c$ in assumption A2. Using compactness of $V^{-1}([0,c])$ and continuity of $V$,  it is easy to show that for any $\epsilon>0$ there is a $b = b(\epsilon)>0$ such that  $V^{-1}([0,b]) \subset S(\Lambda, \epsilon)$. It is also easy to show that for any $b>0$ there is an $\epsilon = \epsilon (b)>0$ such that $S(\Lambda, \epsilon) \subset V^{-1}([0,b])$.  Hence, we can use the family of  sets $\{ V^{-1}([0,b]), b>0\}$ instead of the family of open sets $\{ S(\Lambda, \epsilon), \epsilon >0\}$ to show uniform stability and uniform asymptotic stability of $\Lambda$ for (\ref{our:new:dyn:sys}).

Let us first show uniform stability of $\Lambda$ for (\ref{our:new:dyn:sys}). Given any $\epsilon>0$, take any  $\delta$ such that $0< \delta \leq \min\{\epsilon, c\}$. Then, for any $x_0 \in V^{-1}([0,\delta])$, $x(t; x_0) \in  V^{-1}([0,\delta]) \subset V^{-1}([0,\epsilon])$ for all $t\geq0$ since $V$ is decreasing along the trajectory of $x(t;x_0)$ of  (\ref{our:new:dyn:sys}). Hence, $\Lambda$ is uniformly stable for (\ref{our:new:dyn:sys}). 

Let us now show uniform asymptotic stability of $\Lambda$ for (\ref{our:new:dyn:sys}). Take any $\delta_0$  such that  $0 < \delta_0 \leq c$. By continuous dependence of $x(t;x_0)$ on  initial point $x_0$,  compactness of $V^{-1}([0,\delta_0])$, continuity of the function $V$, and the property that $V(x(t;x_0))$ decreases to 0 as $t\rightarrow \infty$ for any $x_0 \in V^{-1}([0,c])$, it is easy to show that for any $\epsilon>0$ there is a time $T(\epsilon)>0$ such that for any $x_0 \in V^{-1}([0,\delta_0])$ we have $x(t;x_0) \in V^{-1}([0,\epsilon])$ for all $t\geq T(\epsilon)$. Hence, $\Lambda$ is uniformly asymptotically stable for  (\ref{our:new:dyn:sys}). 

\end{proof}
\end{lemma}

Suppose the vector field $X$ is $C^p$ and the function $V$ is $C^{p+1}$ in the modified dynamical system (\ref{our:new:dyn:sys}).  Consider a discrete analogue of  (\ref{our:new:dyn:sys}) described by any one-step numerical method of $p$th order 
\begin{equation}\label{discrete:system}
x_{k+1} = x_k + hY_h(x_k)
\end{equation}
with uniform step size $h>0$, where  $Y_h: \mathbb R^n \rightarrow \mathbb R^n$ for each $h$.
\begin{theorem}
Suppose that the vector field $X$ is $C^p$ and the function $V$ is $C^{p+1}$, and that assumptions A1, A2 and  A3 (or  A3$^{\prime}$ instead of A3) are satisfied. Then there is a number $h_2>0$ such that for each $0<h<h_2$ the discrete-time dynamical system (\ref{discrete:system})  has a compact, uniformly asymptotically stable set $\Lambda_h$ which contains $V^{-1}(0)$ in its interior and converges to $V^{-1}(0)$ with respect to the Hausdorff metric as $h\rightarrow 0+$.  Moreover, there is a bounded, open set $U_0$, which is independent of $h$ and contains $\Lambda_h$, and a time 
\[
T_0(h) = A + B p \log \frac{1}{h},
\]
where $A$ and $B$ are constants depending on the stability characteristic of $V^{-1}(0)$, such that the iterates of (\ref{discrete:system}) satisfy
\[
x_k \in \Lambda_h
\]
for all $kh \geq T_0(h)$, $x_0 \in U_0$ and $0 < h<h_2$.
\begin{proof}
We have only to show that the hypotheses in Theorem 1.1 of \cite{KlLo86} hold. Since $X$ is $C^p$ and $V$ is $C^{p+1}$, the vector field $X - \nabla V$ of  (\ref{our:new:dyn:sys}) and its derivatives of order up to $p$ are all continuous and bounded on the compact set $V^{-1}([0,c])$.  The set $V^{-1}(0)$ is uniformly asymptotically stable for (\ref{our:new:dyn:sys})  by Lemma \ref{lemma:UAS:Lambda} in the above. Therefore, the conclusions of this theorem follow from Theorem 1.1 and Lemma 3.3 of  \cite{KlLo86}.
\end{proof}
\end{theorem}
Refer to \cite{KlLo86} to see how to obtain the set $U_0$ and values of the parameters $h_2$, $A$ and $B$ that appear in the statement of the above theorem. 
The above theorem extends to multi-step numerical integrators; refer to \cite{KlLo90} for detail.

\section*{Acknowledgement}
This research was supported in part by DGIST Research and Development Program (CPS Global Center) funded by the Ministry of Science, ICT \& Future Planning, Global Research Laboratory Program (2013K1A1A2A02078326) through NRF, and Institute for Information \& Communications Technology Promotion (IITP) grant funded by the Korean government (MSIP) (No. B0101-15-0557, Resilient Cyber-Physical Systems Research).

%
%

%
%

%
%

\end{document}